\documentclass[amsthm]{amsart}
\usepackage{graphics}
\usepackage{amsmath}
\usepackage{amscd}
\usepackage{graphicx}
\usepackage{latexsym}
\usepackage{caption}
\usepackage{subcaption}
\usepackage{algorithmicx}
\usepackage{algpseudocode}
\usepackage{algorithm}

\newtheorem{thm}{Theorem}[section]
\newtheorem{Lem}[thm]{Lemma}

\newtheorem{defn}[thm]{Definition}

\newtheorem{Rm}[thm]{Remark}
\newtheorem{question}[thm]{Question}

\def\R{{\mathbf R}}

\def\zr{{\mathbf R}}

\def\citep{\cite}

\def\s2x{\hbox{$S^2 \times S^2$}}

\newcommand{\dvf}{discrete vector field}

\newcommand{\dgvf}{discrete gradient vector field}

\newcommand{\thing}{cell complex}

\newcommand{\lra}{\longrightarrow}

\begin{document}


\title{Birth and death in discrete Morse theory}

\thanks{Third author partially supported by the Research
Agency of Slovenia.}

\author{Henry King}
\address {Department of Mathematics, University of Maryland, College Park,  MD 20742 }
\email{hck@math.umd.edu}

\author{Kevin Knudson} 
\address{Department  of Mathematics, University of Florida, Gainesville, FL 32611}
\email{kknudson@ufl.edu}

\author{Ne\v za Mramor Kosta} 
\address{Department  of Computer and Information Science and Institute of Mathematics, Physics, and Mechanics,
University of Ljubljana,  Slovenia}
\email{neza.mramor@fri.uni-lj.si }

\begin{abstract}
Suppose $M$ is a finite cell decomposition of a space $X$ and that for
$0=t_0<t_1<\cdots <t_r=1$ we have a discrete Morse function
$F_{t_i}:M\to \zr$.  In this paper, we study the births and deaths
of critical cells for the functions $F_{t_i}$ and present an
algorithm for pairing the cells that occur in adjacent slices.  We
first study the case where the cell decomposition of $X$ is the same
for each $t_i$, and then generalize to the case where they
may differ.  This has potential applications in
topological data analysis, where one has function values at a sample of points
in some region in space at several different times or at different
levels in an object.
\end{abstract}



\maketitle

\section{Introduction}
The purpose of this paper is to study the discrete analogue of the
following phenomenon in classical smooth Morse theory.  Suppose
that $N$ is a smooth manifold and that we have a family of
functions $f:N\times I\to\zr$ such that the various $f_t:N\to\zr$
are generically Morse; that is, for almost all $t$, the function
$f_t$ has only nondegenerate critical points.  Then as $t$ varies,
the critical points of the $f_t$ move around in $N$.  Sometimes a
critical point is ``born"; that is, a new critical point appears
at some time $t_0$.  At other times, critical points ``die".
Generically, critical points are born and die in pairs. Such
events are isolated since the critical points of a Morse function
are separated; we call such points in $N\times I$ {\em birth-death
points}.

\begin{figure}
\centerline{\includegraphics[width=4.5cm]{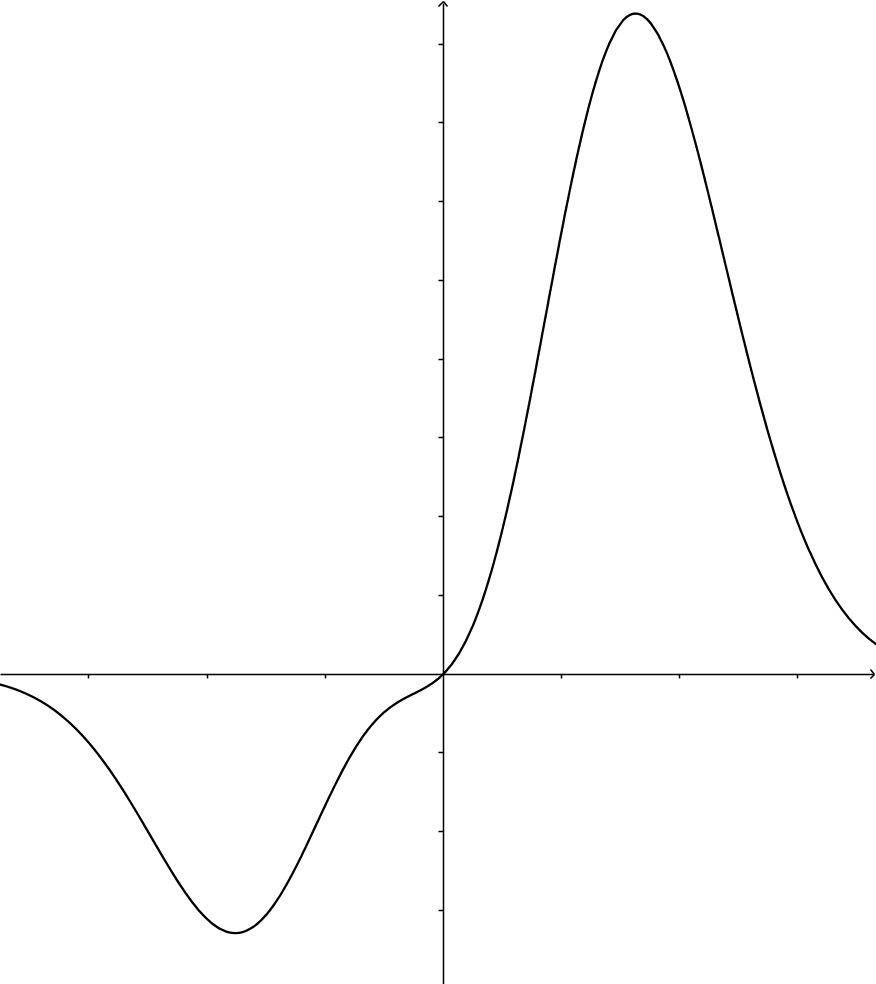}\hskip 3mm
\includegraphics[width=4.5cm]{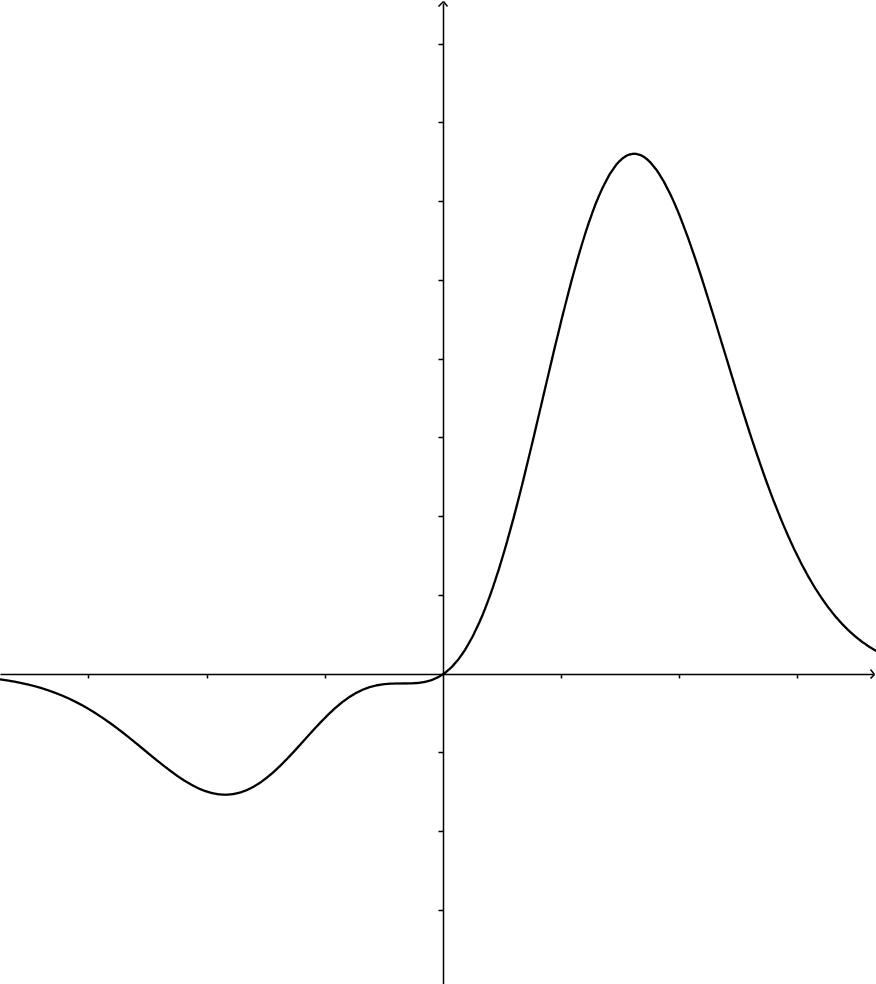} \hskip 3mm
\includegraphics[width=4.5cm]{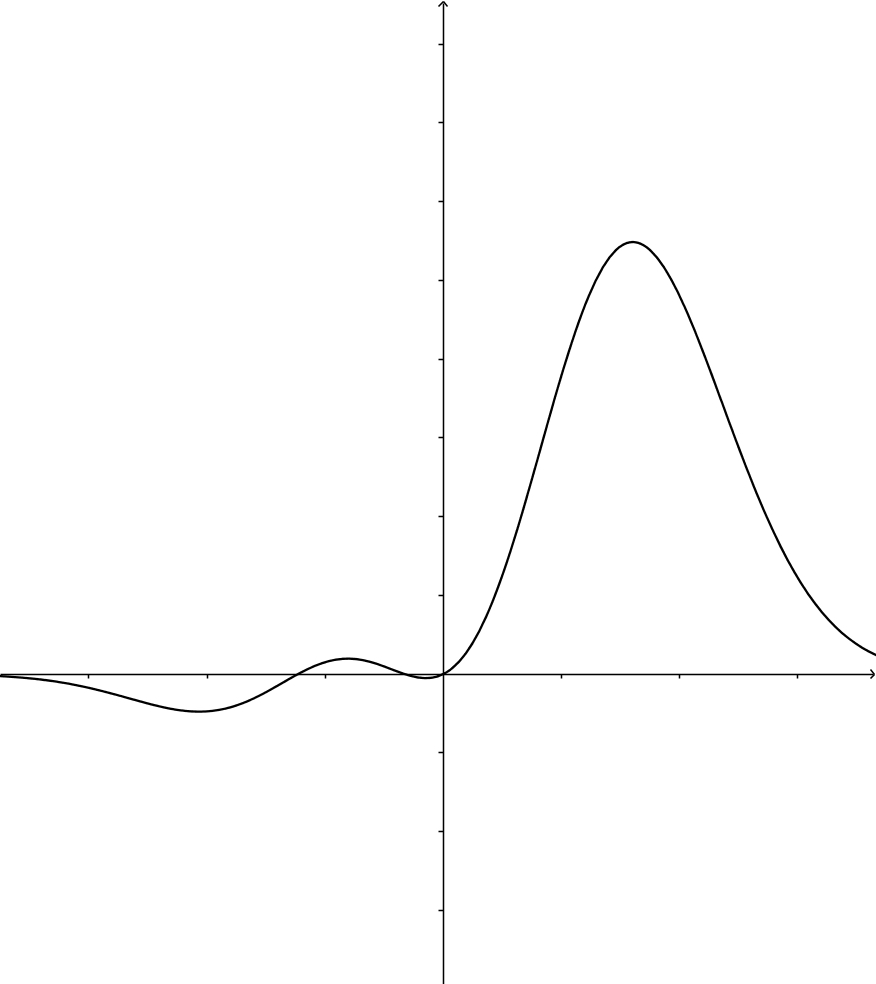}}
\vskip 3mm 
\centerline{\includegraphics[width=4.5cm]{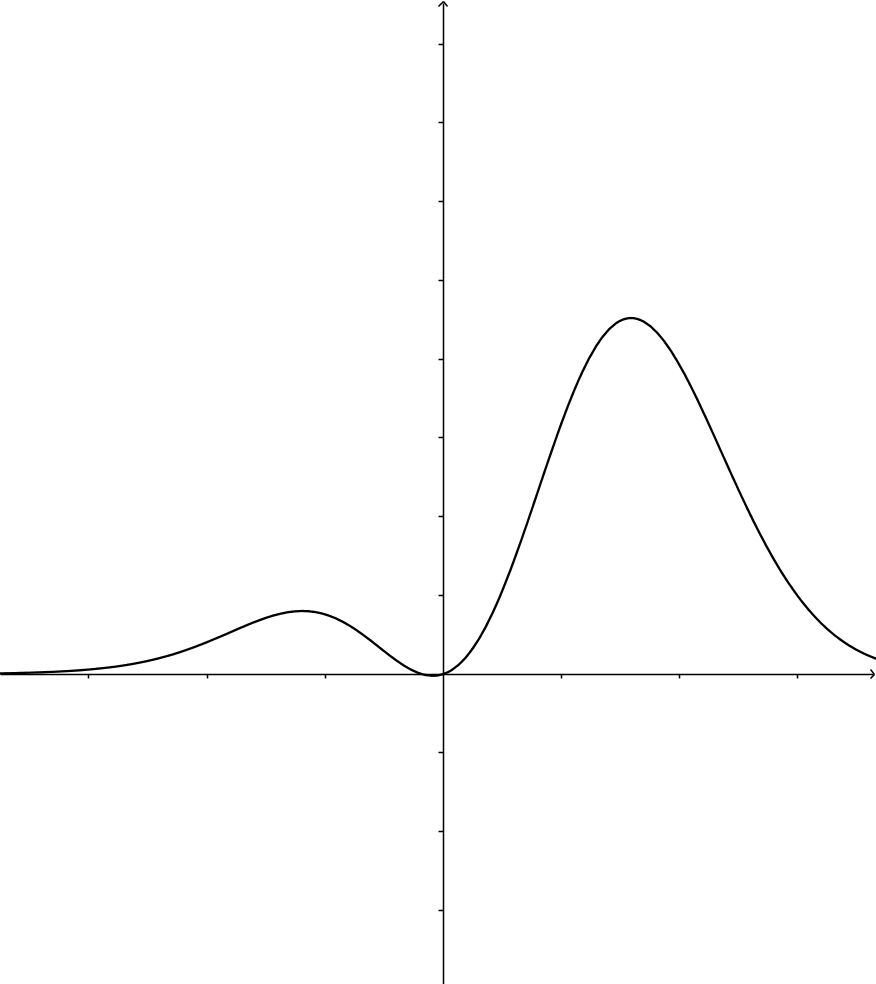}\hskip 3mm
\includegraphics[width=4.5cm]{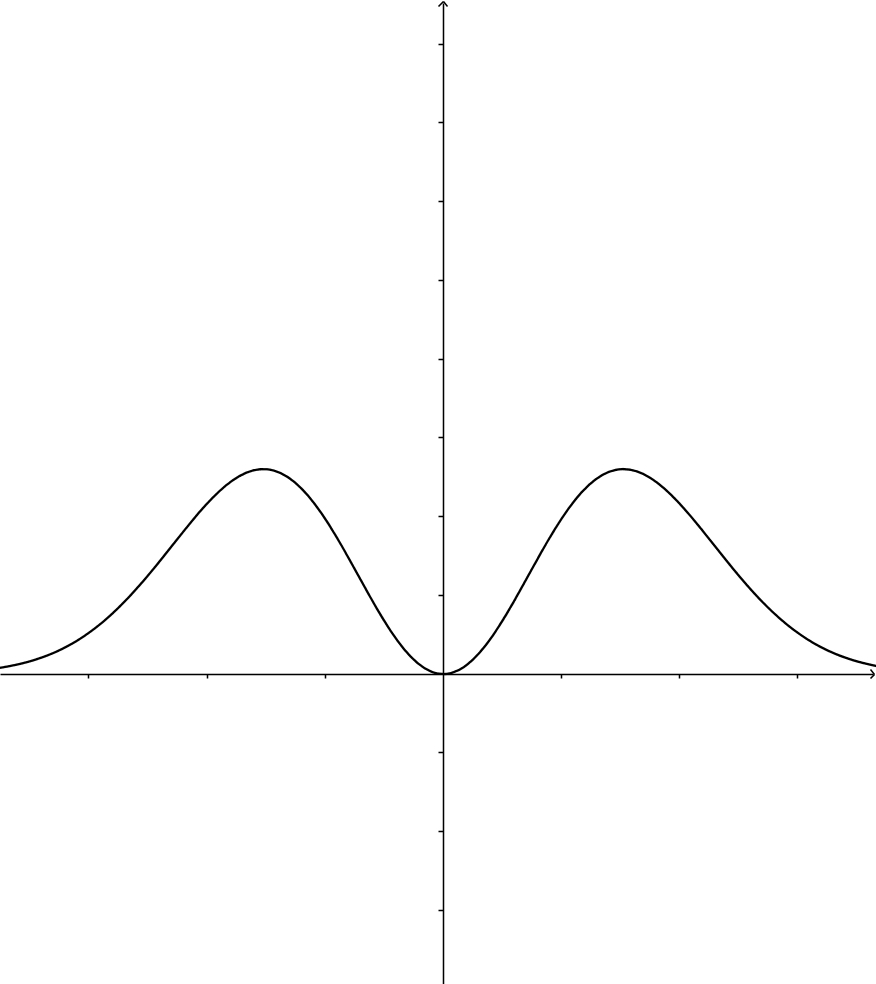}\hskip 3mm
\includegraphics[width=4.5cm]{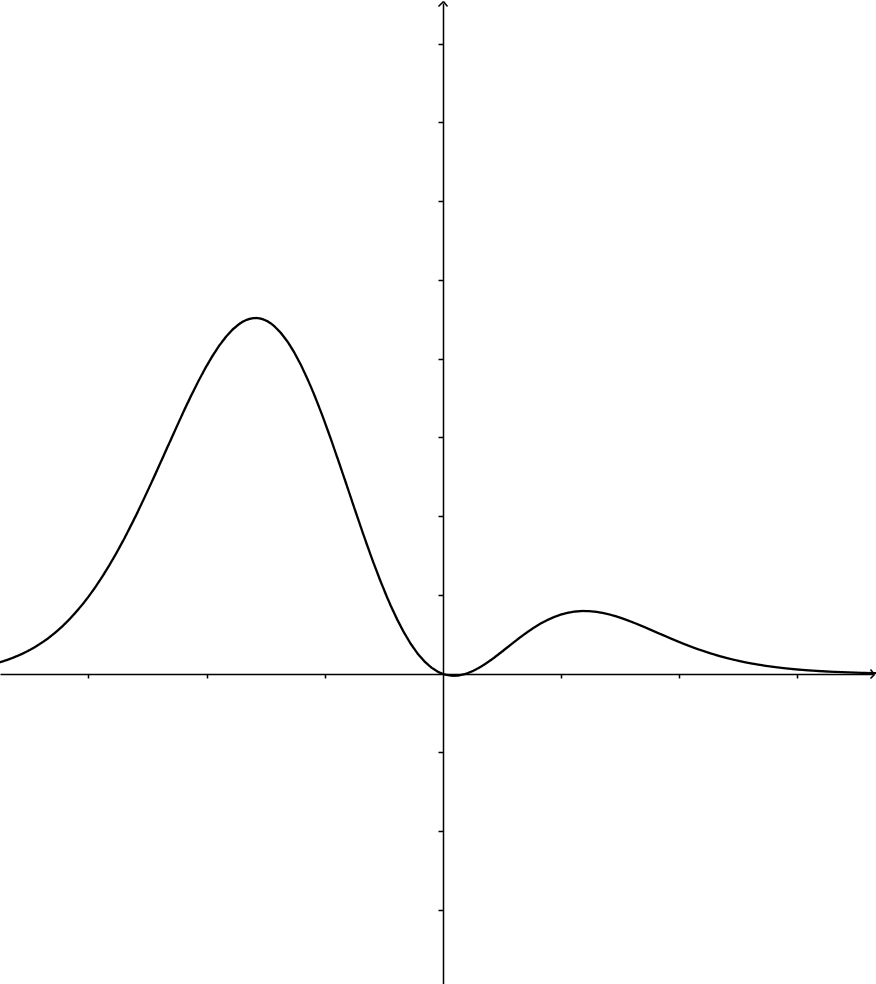}}\hskip 3mm
\vskip 3mm
\centerline{\includegraphics[width=4.5cm]{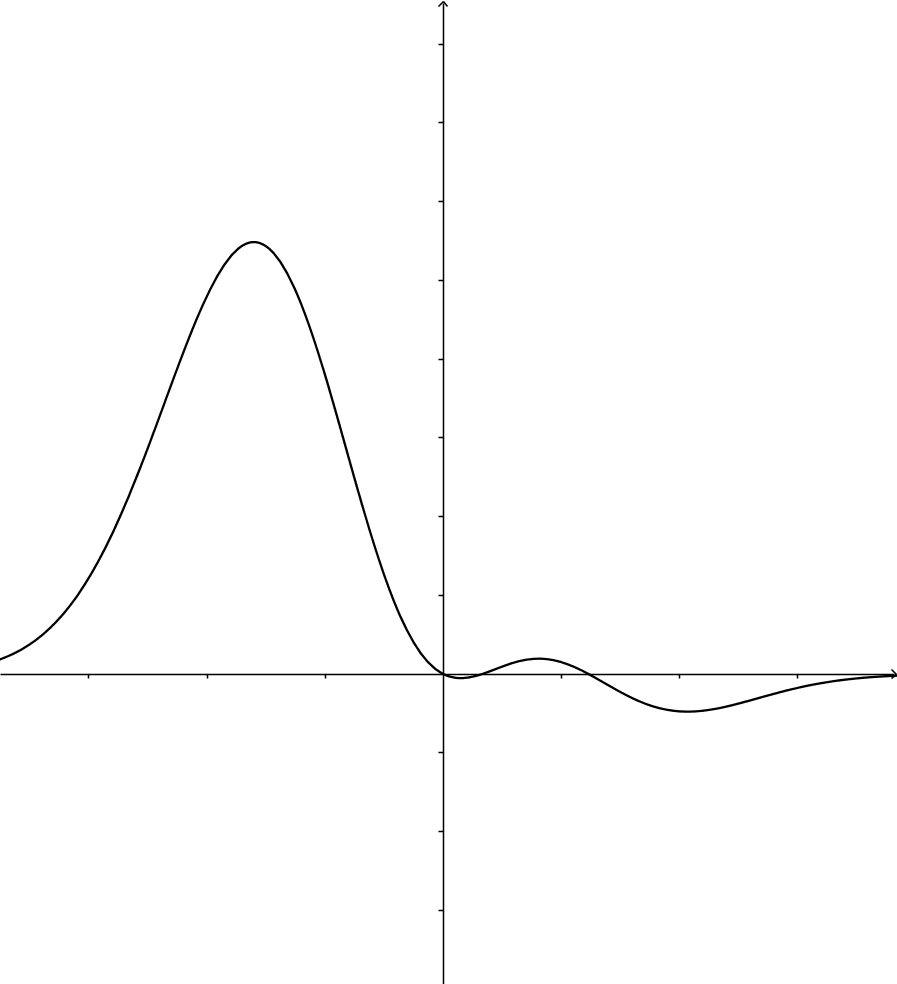}\hskip 3mm
\includegraphics[width=4.5cm]{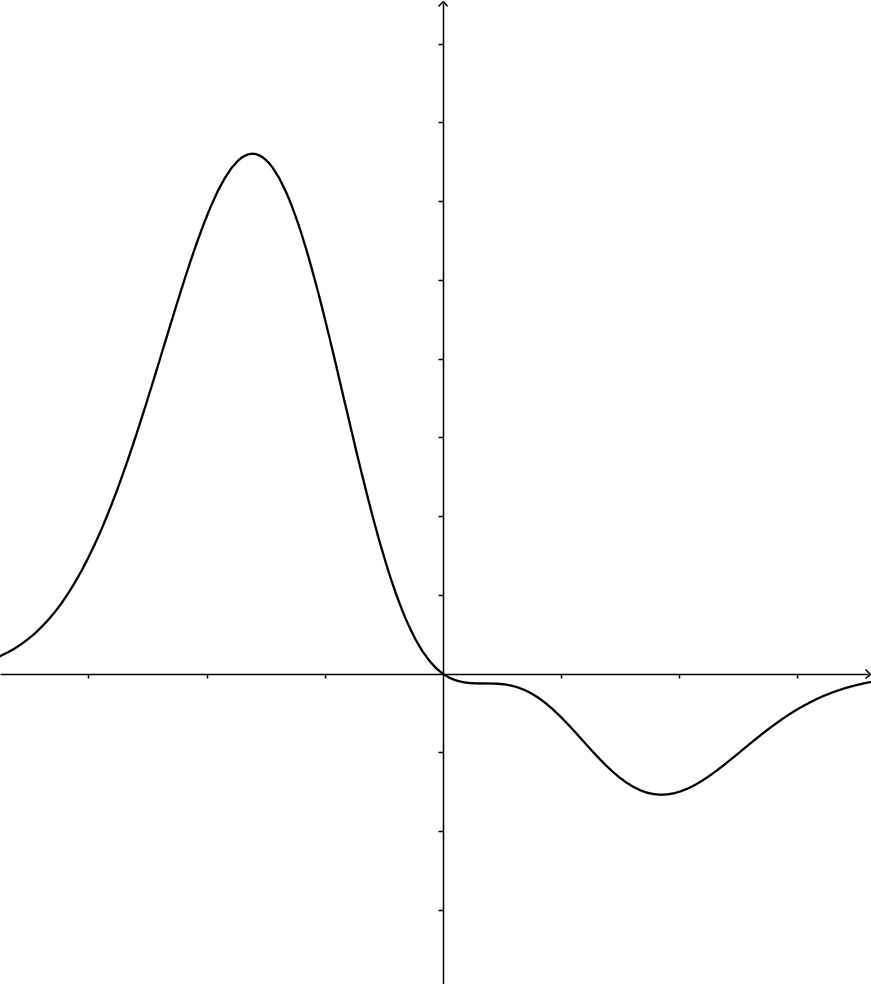}\hskip 3mm
\includegraphics[width=4.5cm]{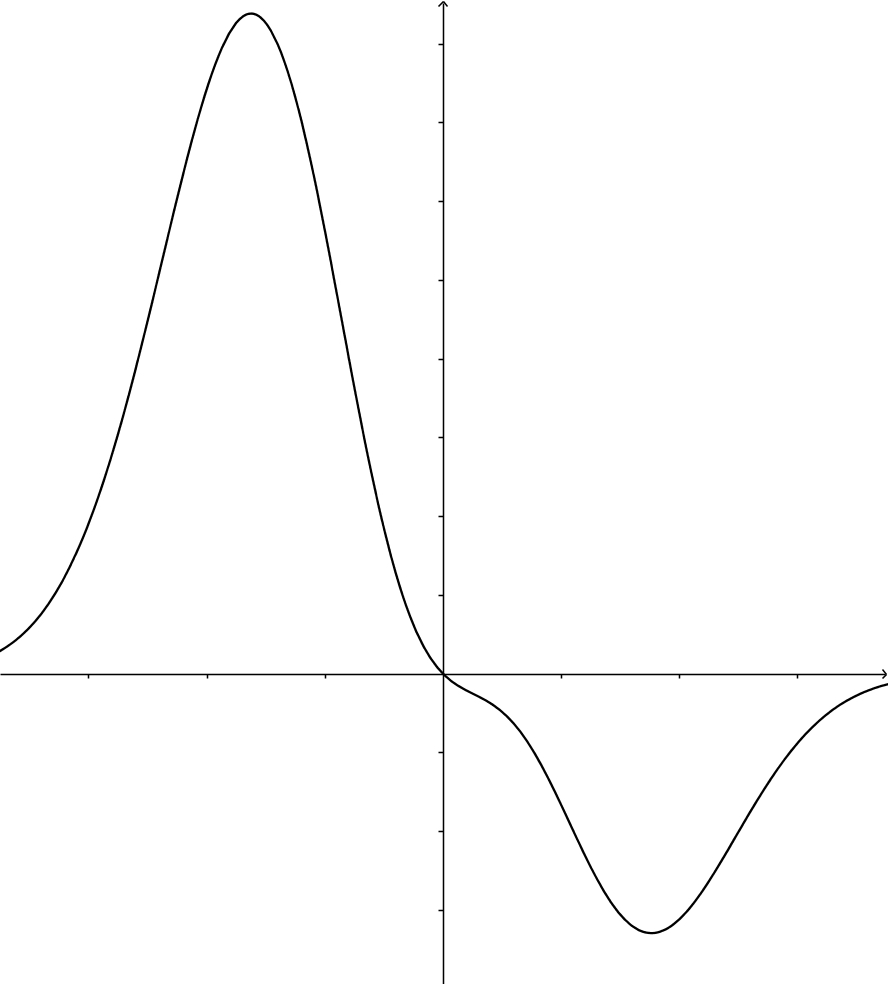}}
\caption{\label{curves} The curves $f(x,t)$ for $t=-3,-2.08,-1.5,-1,0,1,1.5,2.08,3$.}
\end{figure}

Consider the following simple example.  Let
$N=\zr$ and consider the family $f:N\times I\to\zr$ defined
by
$$f(x,t)=e^{-x^2/2}\biggl(\frac{x^4}{2}-3tx^3+6x^2-tx\biggr).$$ Each $f_t$ is a
Morse function except for
a pair of values $t\approx \pm 2.08$.
Figure \ref{curves} shows the graphs of $f(x,t)$ for a
few values of $t$.  Note the evolution of the critical points as we pass through various $t$ values.  In intermediate stages, we see the appearance of degenerate critical points ($t\approx -2.08$ and $t\approx 2.08$).

\begin{figure}
\centerline{\includegraphics[width=4.0in]{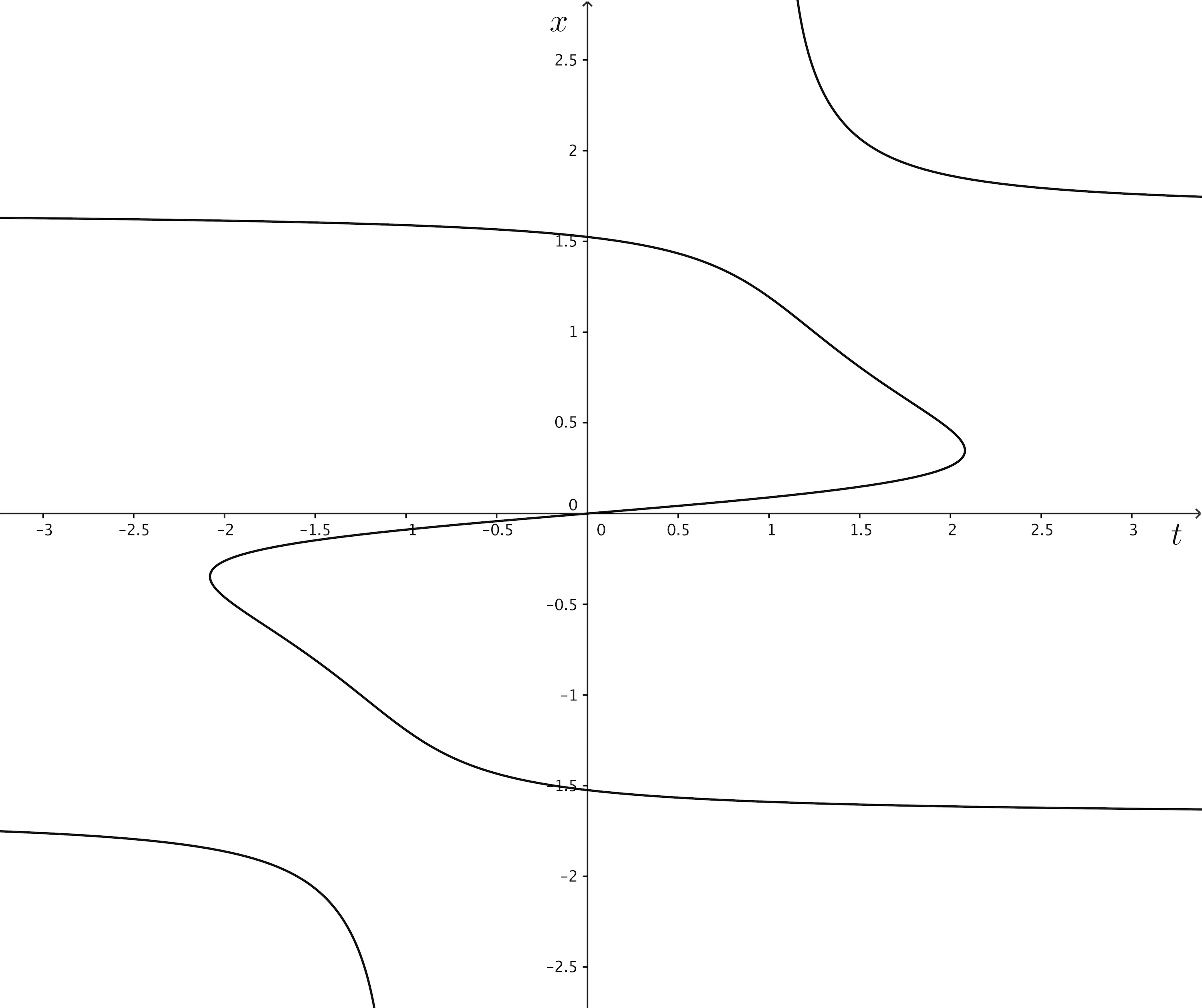}}
\caption{\label{curvebif} The bifurcation diagram associated to
$e^{-x^2/2}(x^4/2-3tx^3+6x^2-tx)$; $t$ is the horizontal axis, $x$ is the vertical.}
\end{figure}

We can trace the critical points as $t$ varies via a bifurcation
diagram as in Figure \ref{curvebif}.  Note that in this simple example, the
diagram is merely the graph of $f_x(x,t)=0$ with $x$ on the vertical
axis and $t$ on the horizontal.

We are concerned here with the following question.  Suppose we
have a cell decomposition $M$ of a space $X$ and for each $0=t_0<t_1<\cdots
<t_r=1$ we have a discrete Morse function $F_{t_i}:M\to\zr$ (the
basics of discrete Morse theory are given in Section \ref{dmt}). We wish to study how
the critical cells move in $M$ as $i$ varies and construct
bifurcation diagrams similar to that in Figure 1.2. To do this, we
need to  develop an algorithm to
pair critical cells in each slice with those in the next one.
We carry this out first in the case where the cell decomposition $M$ of $X$
is the same for each $F_{t_i}$ (Section \ref{alg1}) and then
extend this to the more general case where we allow the possibility that we have chosen a different
cell decomposition $M_i$ of $X$ for each $i=0,\dots ,r$, assuming that for each $i$, $M_i$ and $M_{i+1}$ have a common subdivision (Section \ref{alg2}).  Section \ref{things1} contains general results about discrete gradient vector fields and refinements of them on subdivisions of complexes; these may be of independent interest to practitioners.  In Section \ref{smooth} we justify our algorithms by explaining the parallels with families of smooth functions.  In Section \ref{apps}, we show the results of running our algorithms
on CT scans of the head, where the function in question is density.
Minima correspond to cavities in the interior of the head.  Also, in
a series of renal scintigrams we trace the maxima of the intensity
of radiation which correspond to the locations of the kidneys and
bladder within the image.

\section{Discrete Morse Theory}\label{dmt}

In this section we review the
basics of Forman's discrete Morse theory \citep{F}.
In topological data analysis the domain on which the data is given is often
a simplicial complex, obtained using some triangulation algorithm on the given
data points. It will be convenient in applications as well as in proofs
to consider more general cell complexes, though.

Let $M$ be a cell complex, which for us will mean a
regular CW complex,
 i.e., a CW complex where each cell attaching map is an embedding onto a subcomplex.
A typical $p$-cell will be denoted by
$\alpha^{(p)}$, and if $\alpha$ is a codimension-one face of $\tau$, we write $\alpha<\tau$.
In the following definition, we shall abuse notation slightly
by speaking of functions from $M$ to $\R$ when we really mean functions from
the set of cells of $M$ to $\R$.  No confusion should result.

\begin{defn} A function $f:M\lra\R$ is a {\em discrete Morse function} if
for every $\alpha^{(p)}\in M$, the following two conditions hold:
\begin{enumerate}
\item $\#\{\beta^{(p+1)}>\alpha| f(\beta)\le f(\alpha)\}\le 1,$
\item $\#\{\gamma^{(p-1)}<\alpha| f(\gamma)\ge f(\alpha)\}\le 1.$
\end{enumerate}
\end{defn}

Essentially, then, discrete Morse functions are functions on $M$
that increase with the dimension of the cells; that is, the
values on all but at most one face of $\alpha$ must be smaller
than the value on $\alpha$ itself.  Discrete Morse functions
exist.  Indeed, the simplest (and most trivial) example is the
following.  If $\sigma$ is a cell of $M$, we define $f:M\lra
\zr$ by
$$f(\sigma)=\dim\sigma.$$
Note that the conditions (1) and (2) in the definition of a
discrete Morse function are exclusive.
That is, if one of the sets has cardinality one, the other is empty.
To see this, assume that $\sigma^{(p)}
\in M$ has a coface $\tau^{(p+1)}$ such that $f(\tau)\leq
f(\sigma)$ and a face $\upsilon^{(p-1)}$ such that
$f(\upsilon)\geq f(\sigma)$. Let $\sigma'^{(p)}$ be a different
face of $\tau$ which also contains $\upsilon$. Then the value
$f(\sigma')$ must be bigger than $f(\upsilon)$, since $\upsilon$
already has a coface with smaller value, and similarly the value
of $\tau$ must be bigger than $f(\sigma')$, so
$$f(\tau)\leq f(\sigma)\leq f(\upsilon)<f(\sigma')<f(\tau)$$
which is not possible.

Smooth Morse theory depends heavily on the notion of a critical
point.  The discrete theory does as well.

\begin{defn} Let $f:M\lra \R$ be a discrete Morse function.  A
cell $\alpha^{(p)}$ is {\em critical} if the following two
conditions hold:
\begin{enumerate}
\item $\#\{\beta^{(p+1)}>\alpha| f(\beta)\le f(\alpha)\}=0,$

\item $\#\{\gamma^{(p-1)}<\alpha| f(\gamma)\ge f(\alpha)\}=0.$
\end{enumerate}
A cell that is not critical is called {\em regular}.
\end{defn}

For example, if $f:M\lra \R$ is given by $f(\sigma)=\dim\sigma$,
then every cell is critical.

\begin{thm}[\citep{F}, Theorem 2.5]\label{collapse} Suppose $M$ is a \thing\ with a discrete Morse function.  Then $M$ is homotopy
equivalent to a CW-complex with exactly one cell of dimension $p$
for each critical cell of dimension $p$.
\end{thm}

Constructing discrete Morse functions with not too many critical cells
is rather involved.
Recall that in smooth Morse theory, one often works with the
gradient vector field of a Morse function.  This encapsulates the
qualitative behavior of the function so that one need not even
know the values of the function itself.  

\begin{figure}
\centerline{\includegraphics[width=2.5in]{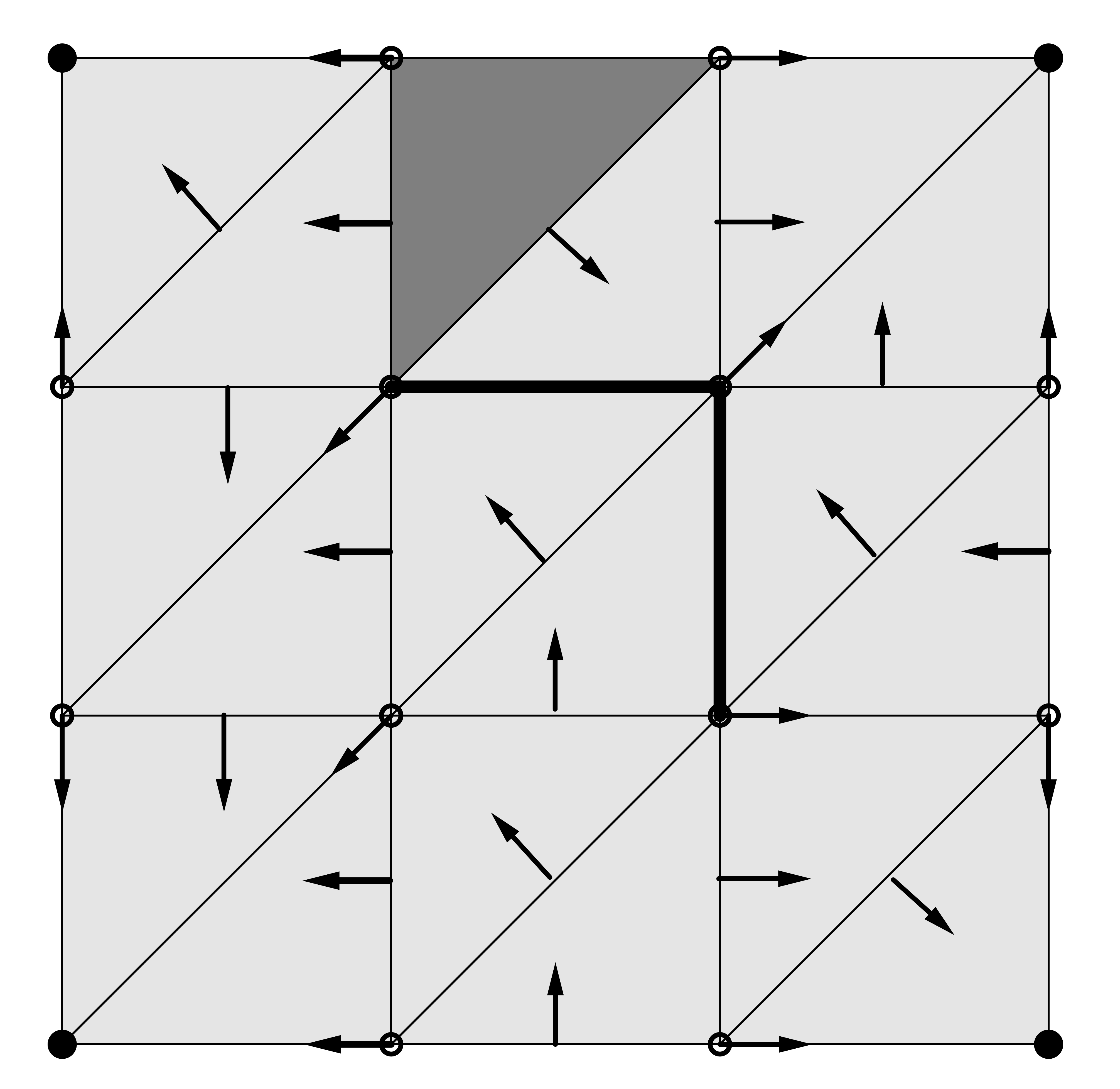}}
\caption{\label{torusfield} A discrete vector field on a triangulated torus.  There are four critical cells: the top-center triangle (dark gray), the top and right edges of the center square (indicated with thicker lines), and the vertex obtained by identifying the four corners.}
\end{figure}

The fundamental observation in the discrete setting is this:  regular cells occur in
pairs.
For an arbitrary cell complex with a discrete Morse function
$f$, we draw arrows to represent a vector field as follows.  If $\alpha^{(p)}$ is a
regular cell with $\beta^{(p+1)}>\alpha$ satisfying
$f(\beta)\le f(\alpha)$, then we draw an arrow from $\alpha$ to
$\beta$.  Figure \ref{torusfield} shows an example of a vector field on the torus (identify the
top side with the bottom, and the left side with the right).

It is easy to see that every cell $\alpha$ satisfies exactly
one of the following:
\begin{enumerate}
\item $\alpha$ is the tail of exactly one arrow;

\item $\alpha$ is the head of exactly one arrow;

\item $\alpha$ is neither the head nor the tail of an arrow.
\end{enumerate}
A cell is critical if and only if it satisfies condition (3)
above.  These arrows can be thought of as the gradient vector
field of the Morse function.  A better point of view is the
following.  A discrete vector field $V$ can be thought of as a
collection of pairs $\{\alpha^{(p)},\beta^{(p+1)}\}$ of cells,
where a pair $\alpha,\beta$ is in $V$ if and only if
$\alpha<\beta$ and $f(\beta)\le f(\alpha)$.

\begin{defn} A {\em discrete vector field $V$} on $M$ is a collection of
pairs $(\alpha^{(p)},\beta^{(p+1)})$ of cells of $M$ with
$\alpha<\beta$ such that each cell is in at most one pair of
$V$.  We often use the notation $V(\alpha)=\beta$ to indicate that the pair $(\alpha,\beta)\in V$.
\end{defn}

Of course, one would like to know if a discrete vector field $V$
is the gradient field of some discrete Morse function on $M$. The
following definition is crucial.

\begin{defn} Let $V$ be a discrete vector field on $M$.  A {\em
V-path} is a sequence of cells
$$\alpha_0^{(p)},\beta_0^{(p+1)},\alpha_1^{(p)},\beta_1^{(p+1)},
\dots ,\beta_r^{(p+1)},\alpha_{r+1}^{(p)}$$ such that for each
$i=0,\dots ,r$, $(\alpha_r,\beta_r)\in V$ and
$\beta_{i}>\alpha_{i+1}\ne\alpha_i$.  Such a path is a {\em
non-trivial closed path} if $r\ge 1$ and $\alpha_0=\alpha_{r+1}$.
The first and last cells are optional, so we also allow a $V$-path
to start or end with a $p+1$ cell.
\end{defn}

For notational convenience in a proof we may assume a path
starts and ends with a lower dimensional $p$-cell, leaving
the other three cases to the reader without notice.

The idea is that a path corresponds to a union of nearby trajectories
of a vector field. Forman proved the following two theorems.

\begin{thm}[\citep{F},Theorem 3.4]  Suppose $V$ is the gradient
vector field of a discrete Morse function $f$.  Then a sequence of
cells is a $V$-path if and only if
$\alpha_i<\beta_i>\alpha_{i+1}$ for $i=0,\dots ,r$, and
$$f(\alpha_0)\ge f(\beta_0)>f(\alpha_1)\ge f(\beta_1)>\cdots \ge
f(\beta_r)>f(\alpha_{r+1}).$$
\end{thm}

In particular, if $V$ is a gradient vector field, then there are
no nontrivial closed $V$-paths.  The converse is true as well.

\begin{thm}[\citep{F},Theorem 3.5]\label{thm:dmfexists} A discrete vector field $V$ is
the gradient vector field of a discrete Morse function if and only
if there are no nontrivial closed $V$-paths.
\end{thm}

\begin{Rm} Forman proved Theorem \ref{thm:dmfexists} only for finite cell complexes.  It is also true for infinite cell complexes; see \citep{infinite1}.
\end{Rm}

The reader is invited to check that the discrete vector field shown
in Figure 2.1 is the gradient of a discrete Morse function on the torus.

A convenient combinatorial description of vector fields may be
given in terms of the Hasse diagram of $M$.  This is the directed
graph whose vertices are the cells of $M$, and whose edges are
given by the face relations in $M$ (i.e., there is an edge from
$\beta$ to $\alpha$ if and only if $\alpha$ is a codimension-one
face of $\beta$). Given a vector field $V$, modify the Hasse
diagram in the following manner.  If $(\alpha,\beta)\in V$, then
reverse the orientation of the edge between $\alpha$ and $\beta$.
A $V$-path is then a directed path in the modified graph.

\begin{thm}[\citep{F},Theorem 6.2] There are no nontrivial closed $V$-paths if and only
if there are no nontrivial closed directed paths in the modified
Hasse diagram.
\end{thm}

Given a finite cell complex $M$ and a real valued function on its $0$-cells,
the algorithm in \citep{kkm} gives a \dgvf\ which is meant to
be a combinatorial analogue of the gradient vector field
of a function sampled at the $0$-cells. (Though explicitly stated for simplicial complexes, it easily
extends to cell complexes.)  It may produce additional critical cells that do not correspond to changes in the homotopy type of the space, but this is often unavoidable.  Indeed, producing a gradient with no extraneous critical cells is a hard problem \citep{joswig}.

\section{The Algorithm, Part I}\label{alg1}
As a special case, in this section we assume the following.
Let $M$ be a finite cell decomposition of a space $X$ and suppose that for each $0=t_0<t_1<\cdots <t_r=1$,
we have a
discrete Morse function $F_{t_i}:M\to\zr$.
(Note we always use the {\em same} cellulation
$M$ for each $i$.)
In applications, we would probably have the values of each $F_{t_i}$
only on the zero cells of $M$, but we can always extend this to a discrete Morse
function on all of $M$ via the algorithm we presented in \citep{kkm}.  There is no canonical choice of such a function, but it may be taken to be arbitrarily close to the function assigning to each cell the maximum of the values among its vertices.  We therefore assume
that this has already been accomplished and let $V_i$
denote the gradient vector field for $F_{t_i}$.
We think of  $F_{t_i}$ as a function which varies with time
and we wish to find which critical cells for $V_i$ correspond
to which critical cells for $V_{i+1}$, similar to the smooth example presented in the introduction.

\begin{defn}\label{strongconnectdef} Suppose $\alpha$ and $\beta$ are $k$-cells of $M$
and $\alpha$ is critical for $V_i$ and $\beta$ is critical for $V_j$,
$j\ne i$.  We say that $\alpha$ is {\em connected} to $\beta$ if there
is a $k$-cell $\gamma$ and a $V_i$-path $\alpha,\ldots,\gamma$
of $k$- and $(k-1)$-cells and a $V_j$-path $\gamma,\ldots,\beta$
of $k$- and $(k+1)$-cells.
We say that $\alpha$ is {\em strongly connected} to $\beta$
if $\alpha$ is connected to $\beta$ and $\beta$ is connected to $\alpha$.
\end{defn}

For example if $k=0$ we must have $\gamma=\alpha$
so a $0$-cell $\alpha$ is connected to $\beta$ if there is a
$V_j$-path from $\alpha$ to $\beta$.
There is necessarily exactly one such path.
However, higher dimensional critical cells could be
connected to more than one cell, or even to no cell.

Algorithm \ref{fc} below returns for each pair $i,i+1$ all pairs of connected critical cells. The input to the algorithm are four lists obtained from the algorithm in \citep{kkm}: the lists $C_i$ and $C_{i+1}$ of critical cells in slices $i$ and $i+1$, and the lists of pairs of $V_i$ and $V_{i+1}$. 

\begin{algorithm}
\caption{Connecting critical cells}
\label{fc}
\begin{algorithmic}[1]
\For{$d=1,\ldots,\dim M$}
\State Set $B$ to an empty list of pairs
\ForAll{$\alpha^d$ in $C_i$}

\State \Return{list $L$ of pairs $(\tau^{d-1},\sigma^d)$ belonging to a path in $V_i$ starting in $\partial\alpha$}

\ForAll{$\sigma^d$ from $L$} 

\State \Return{list $B_0$ of ending cells $\beta^d$ in $C_{i+1}$ of paths in $V_{i+1}$ starting in $\sigma^{d}$}

\EndFor

\State append list of pairs $(\alpha,\beta)$, $\beta\in B_0$ to $B$

\EndFor
\EndFor
\end{algorithmic}
\end{algorithm}

Strongly connected pairs of critical cells are connected pairs $(\alpha,\beta)$ obtained after running Algorithm \ref{fc} forward as well as backward with respect to time. 

The motivation for our definition of strongly connected critical cells is obtained by examining the corresponding operation in the smooth case: to move from a critical point of a smooth function $f$ to one of a nearby smooth function $g$, first flow along the descending disc of one and then along the ascending disc of the other; these intersect at a single point. In discrete Morse theory we do not have the obvious duality that occurs in the smooth case (the index-$i$ critical points of a Morse function $f$ are the index-($n-i$) critical points of $-f$) and so we must seek an alternate formulation of these discs. The combinatorial analogue of the descending disc for a critical $k$-cell $\sigma$ is the union of all $k$-cells $\tau$ reachable by a $k,(k-1)$-gradient path from $\sigma$.  The ascending disc is the corresponding object in a suitably defined dual cell decomposition, the union of all $n-k$ dual cells reachable from the dual of $\sigma$.  These intersect in a discrete set of points.  More information about the proper definition of the dual complex and the discrete descending and ascending regions is available in \citep{jm}.  In Definition \ref{strongconnectdef}, the path from $\alpha$ to $\gamma$ is in the descending region and the path from $\gamma$ to $\beta$ is the dual of a path in the ascending region from the dual of $\beta$ to the dual of $\gamma$.


\begin{defn}\label{borndied} Suppose $\alpha$ is critical for $V_i$. If $\alpha$ is not strongly connected to any critical cell for $V_{i+1}$, we say $\alpha$ {\em dies} at $i$. If $\alpha$ for $V_i$ is not strongly connected
to any critical cell for $V_{i-1}$ we say $\alpha$ is {\em born} at $i$. If $\alpha$ is strongly connected to exactly one critical cell $\beta$ for $V_{i+1}$
and $\beta$ is not strongly connected to any other
critical cell for $V_i$ we  say $\alpha$ {\em moves} to $\beta$.
\end{defn}

Given these connectivity data, we may then draw diagrams illustrating connections among critical cells analogous to the bifurcation diagrams generated by a smooth family of functions.  Algorithm \ref{birthdeathalg1} returns a collection of such diagrams, one for each dimension $d=0,\ldots, \dim  M$. 

\begin{algorithm}[H]
\caption{Algorithm Birth-death \#1}
\label{birthdeathalg1}
\begin{algorithmic}[1]
\For{$d=0,\ldots,\dim M$}
\For{$i=0,\ldots, r$}
\ForAll{$\alpha^d$ critical in $V_i$} 
\State add a node (denoted by $\alpha$) in layer $i$ 
\If{$i>0$}
\ForAll{$\alpha^d$ critical in $V_i$}
\State add an edge connecting the node $\alpha$ to each node in $V_{i-1}$, strongly connected to $\alpha$
\EndFor
\EndIf
\EndFor
\EndFor
\EndFor
\end{algorithmic}
\end{algorithm}

\subsection*{A Simple Example}

\begin{figure}
\centerline{\includegraphics[width=2.5in]{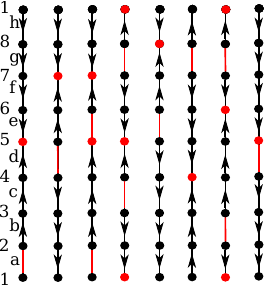}}
\caption{\label{alg1fig1} A family of discrete gradient fields on the circle. The slices are numbered left to right as $0,1,\dots ,7$.}
\end{figure}

\begin{figure}
\centering
\begin{subfigure}{.5\textwidth}
  \centering
  \includegraphics[width=.9\linewidth]{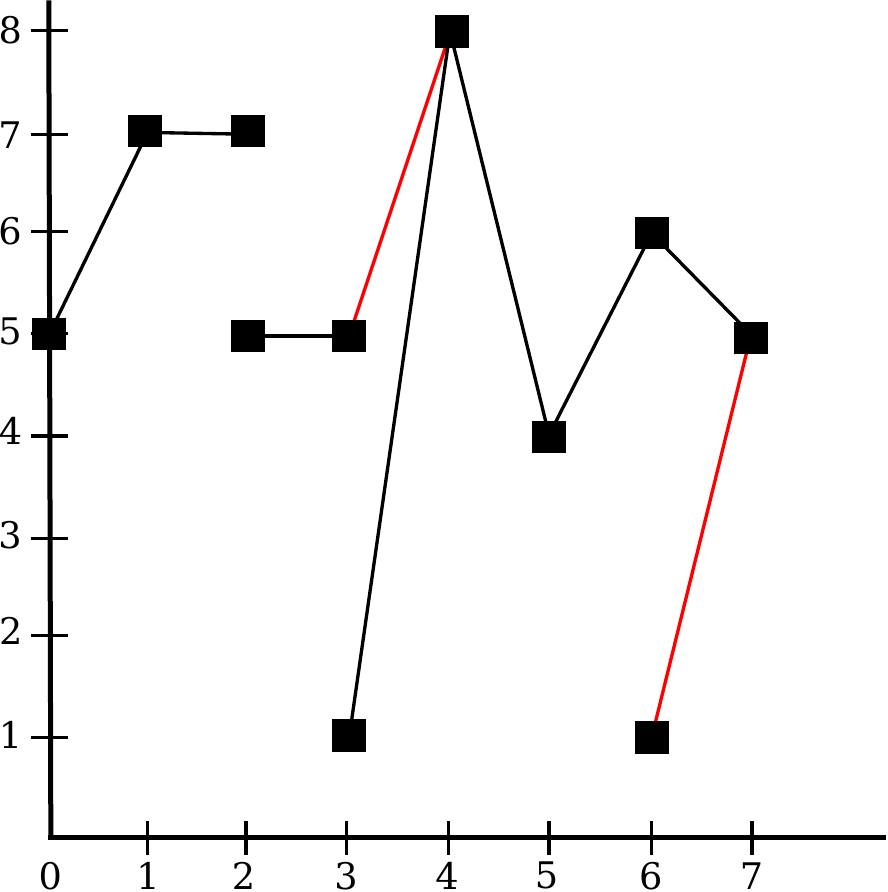}
  \caption{$0$-cells}
\end{subfigure}%
\begin{subfigure}{.5\textwidth}
  \centering
  \includegraphics[width=.9\linewidth]{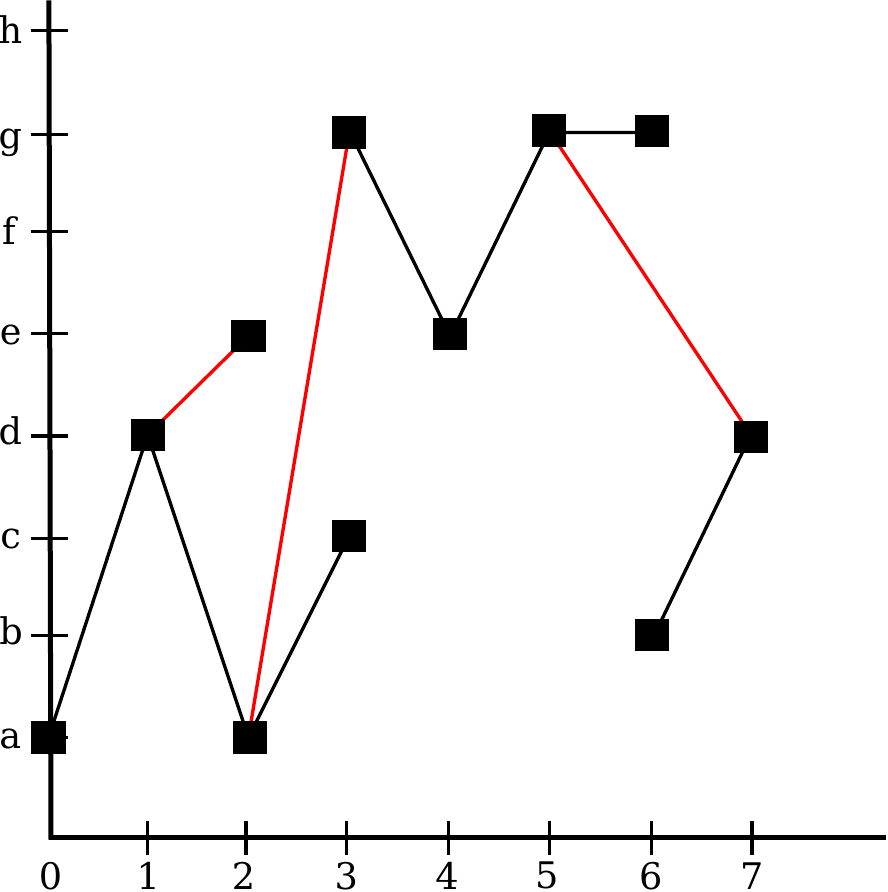}
  \caption{$1$-cells}
\end{subfigure}
\caption{\label{alg1fig2} The connectivity diagrams associated with the example in Figure \ref{alg1fig1}.}

\end{figure}

In Figure \ref{alg1fig1} we present an example of a family of discrete gradients on the circle.  Note that we have used the same triangulation in each case.  Figure \ref{alg1fig2} shows the connectivity diagrams for the $0$- and $1$-cells.  Strong connections are indicated by black lines while connections (in the forward direction only) in adjacent slices are shown in red.  We have the following births:  vertex $5$ in slice $2$, vertex $1$ in slice $3$, and vertex $1$ in slice $6$; edge $e$ in slice $2$, edge $g$ in slice $3$, and edge $b$ in slice $6$.  The deaths are: vertex $7$ in slice $2$, vertex $5$ in slice $3$, and vertex $1$ in slice $6$; edge $d$ in slice $1$, edge $e$ in slice $2$, and edge $c$ in slice $3$.

\begin{Rm} In many of the bifurcation diagrams that follow, we often show connections that occur in one direction (forward or backward). We do this to show that there can be connections in one direction which are not strong connections. Moreover, while Definition \ref{strongconnectdef} does not care about the relative order of $i$ and $j$, the slice parameter is often thought of as time and so we think it is natural to show forward connections. 
\end{Rm}

\section{Refining discrete vector fields}\label{things1}

To deal with a collection $\{F_{t_i}\}$ of discrete Morse functions defined on different cell decompositions $M_i$, $i=0,1,\dots ,r$, we need to restrict the types of cell complexes we work with.  In the remainder of the paper, we assume that each cell decomposition $M_i$ is either a simplicial complex or a cubical complex, where by cubical complex we mean a polyhedral complex consisting of cubes (but not necessarily having the property that the intersection of two cubes is a common face of both).  These complexes suffice for most data analysis applications.


\begin{defn}\label{cellrefine} A cellulation $N$ of a space is a {\em refinement} of a cellulation $M$
if each cell of $N$ is contained in some cell of $M$.  
\end{defn}


\begin{defn}\label{dvfrefinement}
Let $V$ be a \dvf\ on a \thing\ $M$ and let $N$ be a refinement
of $M$.  If $\sigma$ is a cell of $N$ let $g(\sigma)$
denote the smallest dimensional cell of $M$
so that $\sigma\subset g(\sigma)$.
We say that a \dvf\ $W$ on $N$ is a {\em refinement}
of $V$ if for each $k$-cell $\alpha$ of $M$ we can choose a $k$-cell $h(\alpha)$ of $N$ so that:
\begin{enumerate}
\item $h(\alpha)\subset \alpha$.
\item If $\alpha$ is critical then $h(\alpha)$ is critical.
\item If $(\alpha,\beta)\in V$  then $(h(\alpha),h(\beta))\in W$.
\item If $(\kappa,\sigma)\in W$
and $\kappa\ne hg(\kappa)$ then $g(\kappa)=g(\sigma)$.
\end{enumerate}
The vector field $W$ may or may not have more critical cells than $V$.
\end{defn}

\begin{figure}
\centerline{\includegraphics[width=5in]{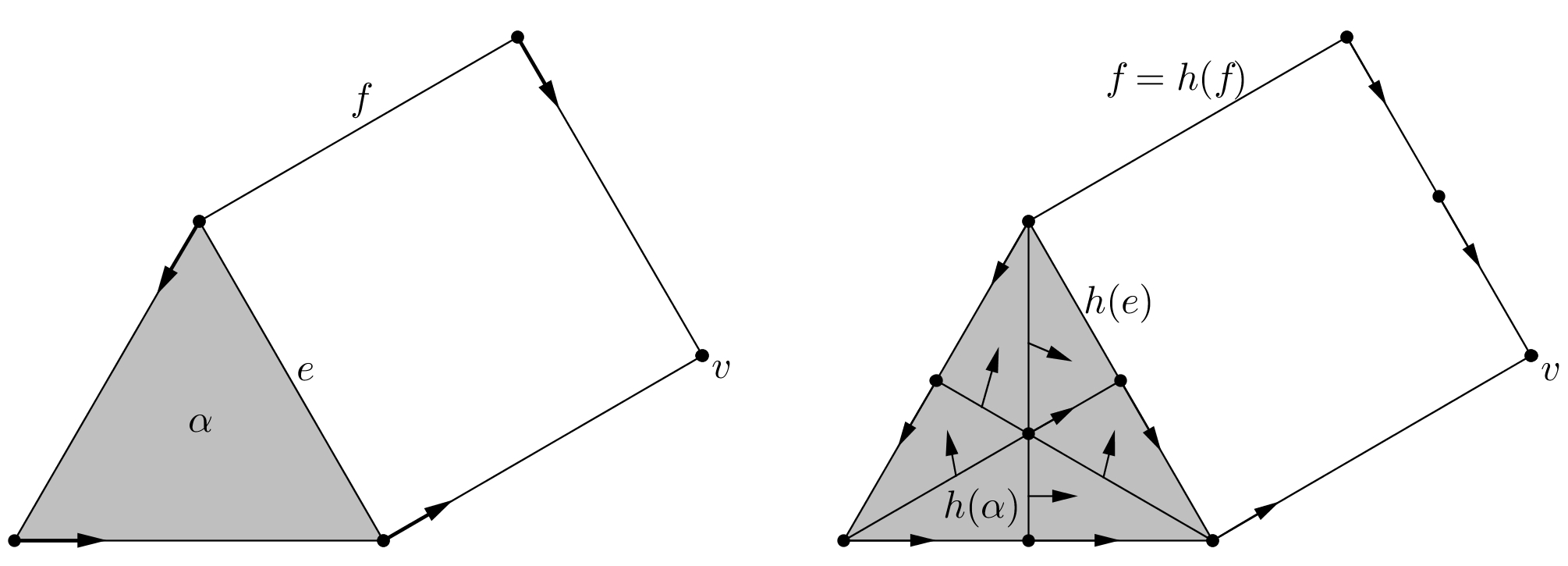}}
\caption{\label{dvrefex} Left: A simplicial complex $M$ with a discrete gradient $V$.  The critical cells of $V$ are $\alpha$, $e$, $f$, and $v$.  Right: A refinement $N$ of $M$ and a refinement $W$ of $V$.  The new critical cells $h(\alpha)$ and $h(e)$ are indicated.}
\end{figure}

An example of a refinement of a discrete vector field is presented in Figure \ref{dvrefex}.  A nonexample could be obtained by taking the midpoint $m$ of the bottom edge of the triangle and pairing it instead with the edge $\tau$ joining $m$ to the barycenter. Then $g(m)$ is the bottom edge $\nu$ in $M$ which has $h(\nu)$ equal to its left half in $N$.  We then have $m\ne hg(m)$ but $g(m) = \nu \ne \alpha = g(\tau)$, violating the fourth condition in Definition \ref{dvfrefinement}. 

Given a \dvf\ $V$ on $M$ and a refinement $N$ of $M$
it is easy to find a refinement $W$ of $V$.
For example, for each $(\alpha^{(p)},\beta^{(p+1)})\in V$,
choose a $p$-cell $\kappa\in N$ with $\kappa\subset \alpha$.
There is a unique $(p+1)$-cell $\sigma$ of $N$ such that
$\sigma\subset \beta$ and $\kappa\subset\sigma$.
Put $(\kappa,\sigma)\in W$.
Such a spare refinement will probably not be useful
because it will contain many critical cells.  We can actually do much better, as we now show, in particular for dimensions at most $2$.

\begin{Lem}
Let $V$ be a \dgvf\ on a  finite \thing\ $M$.
Let $M'$ be a refinement of $M$.
Then there is a \dgvf\ $V'$ on $M'$ which refines $V$.
If $\dim M\le 2$, $V'$ has just one critical cell
for each critical cell of $V$, the minimum possible.
\end{Lem}

\begin{proof}  We suppose by induction on the number of faces of $M$ (where we can proceed dimensionwise) that there is
a subcomplex $K\subset M$ and a function
$h$ from the cells of $K$ to the cells of $M'$
and a \dvf\ $V'$ on $M'$ satisfying
the conditions of the definition of a refinement on $K$.
That is:
\begin{enumerate}
\item $h(\alpha)\subset \alpha$.
\item If $\alpha$ is critical then $h(\alpha)$ is critical.
\item If $(\alpha,\beta)\in V$  and $\beta\subset K$ then $(h(\alpha),h(\beta))\in V'$.
\item If $(\alpha,\beta)\in V$  and $\alpha\subset K$ and $\beta\not\subset K$ then $h(\alpha)$ is critical for $V'$.
\item If $(\kappa,\sigma)\in V'$ then $\sigma\subset K$.
\item If $(\kappa,\sigma)\in V'$
and $\kappa\ne hg(\kappa)$ then $g(\kappa)=g(\sigma)$, where $g$ is the map in Definition 4.3 which assigns to the cell $\tau\in M'$ the smallest cell containing $\tau$ in $M$.
\end{enumerate}
Let $\tau$ be a $k$-cell of $M$ whose boundary is contained
in $K$.
To complete the induction we need to define $h(\tau)$
and also extend $V'$ to cells contained in $\tau$.

First suppose either $\tau$ is critical
or $(\tau,\beta)\in V$ for some $\beta$.
Pick any $k$-cell $\kappa$ of $M'$ in $\tau$
and define $h(\tau)=\kappa$.
This definition of $V'$ satisfies all the inductive conditions,
but it may introduce many extra critical cells.
Instead we apply Lemma \ref{discrefine} below and add $W$ (which is constructed in the proof of Lemma \ref{discrefine}) to $V'$.

Now suppose that $(\beta,\tau)\in V$ for some $\beta$.
Since $\beta$ is in the boundary of $\tau$ we must have
$\beta\subset K$.
We define $h(\tau)$ to be the unique $k$-cell of $M'$
which is contained in $\tau$ and which has $h(\beta)$ in its boundary.
By the inductive assumptions, $h(\beta)$ is critical for $V'$.
Apply Lemma \ref{discrefine} below and add $W$ and $\{(h(\beta),h(\tau))\}$
to $V'$.
\end{proof}

\begin{Lem}\label{discrefine}
Suppose $M$ is a cellulation of the closed $n$-disc
and $\kappa$ is an $n$-cell of $M$.
Then there is a \dgvf\ $W$ on $M$ so that $\kappa$
is critical, all cells in the boundary sphere are critical.
If $n\le 2$ there are no other critical cells.
\end{Lem}

\begin{proof}
The easiest way seems to be a trick.
Consider the dual cell complex $N$ of the interior of $M$.
This is a complex with one $k$-cell $d(\alpha)$
for each $(n-k)$-cell $\alpha$ of $M$ which is not contained
in the boundary sphere.  Also $d(\beta)$ is in the boundary
of $d(\alpha)$ if and only if $\alpha$ is in the boundary of $\beta$.
Choose some distinct values for the $0$-cells of $N$
so that $d(\kappa)$ has the minimum value by far.
Apply the algorithm in \citep{kkm} to get a \dgvf\ $U$
for this function, doing as much cancellation
of critical cells as possible using the algorithm.  

In the end the $0$-cell $d(\kappa)$
will still be critical and if $n\le 2$ it will be the only critical cell.
Now define $W = \{(\alpha,\beta) \mid (d(\beta),d(\alpha))\in U\}$.
Note that $\alpha$ is critical for $W$ if and only if either $d(\alpha)$
is critical for $U$ or $\alpha$ is in the boundary.
Note also that if we take $d$ of a $W$-path and reverse the order, we get a $U$ path; hence $W$ has no nontrivial closed loops.
\end{proof}

\begin{Rm} The cancellation of critical cells mentioned in the proof of Lemma \ref{discrefine} does not maximize the number of possible cancellations, but it cancels in a greedy fashion and in practice yields a \dgvf\ with relatively few critical cells. 
\end{Rm}

\section{The Algorithm, Part II}\label{alg2}
We now consider the general case where we allow different cellulations $M_i$ of the space $X$ for each time slice.
Indeed, this case is important for applications where the
sample of points taken from an object may change over time, thereby yielding
different cell decompositions.  Since these samples form the basis of our algorithm
for generating discrete Morse functions \citep{kkm}, we must consider this
possibility.
However we ask that the cellulations $M_i$ and $M_{i+1}$
have a common subdivision $M_{(2i+1)/2}$.

To see which critical cells for $V_i$
are connected to which critical cells for $V_{i\pm 1}$
we first put a \dvf\ $V^\pm_i$ on $M_{(2i\pm 1)/2}$
which refines $V_i$ and has as few critical cells as
we can manage (see Section \ref{smooth} for a discussion of how well we can do).
Let $h^\pm_i$ and $g^\pm_i$ be the $h$ and $g$
functions appearing in the definition of a \dvf\ refinement.

\begin{defn}\label{strongconnectdef2}
A critical $k$-cell $\alpha$ for $V_i$
is {\em connected} to a critical $k$-cell $\beta$ for $V_{i\pm1}$
if there is a $k$-cell $\gamma$ of $M_{(2i\pm1)/2}$ and
a $V^\pm_i$-path $h^\pm_i(\alpha),\ldots,\gamma$
of $k$- and $(k-1)$-cells and a $V_{i\pm1}$-path $g^\mp_{i\pm 1}(\gamma),\ldots,\beta$
of $k$- and $(k+1)$-cells.  We also want $g^\mp_{i\pm 1}(\gamma)$
to be a $k$-cell.
As before we say that $\alpha$ is {\em strongly connected} to $\beta$
if $\alpha$ is connected to $\beta$ and $\beta$ is connected to $\alpha$.
\end{defn}

Algorithm \ref{fc2} is similar to Algorithm \ref{fc}, with only minor differences. In addition to the four lists  $C_i$, $C_{i+1}$, $V_i$ and $V_{i+1}$ it requires the list of critical cells $C_{(2i+1)/2}$ and regular pairs $V_i^+$ in the common subdivision. 

\begin{algorithm}
\caption{Connecting critical cells between different cell decompositions}
\label{fc2}
\begin{algorithmic}[1]
\For{$k=1,\ldots,\dim M$}
\State Set $B$ to an empty list of pairs
\ForAll{$\alpha^k$ in $C_i$}

\State \Return list $L$ of pairs $(\tau^{k-1},\gamma^k)$ which belong to a path in $V_i^+$ beginning in $\partial h_i^+(\alpha)$ and such that $\dim g^-_{i+1}(\gamma)=k$
\ForAll{$\sigma^k$ from $L$} 

\State \Return list $B_0$ of ending cells $\beta^k$ in $C_{i+1}$ of paths in $V_{i+1}$ beginning in $g^-_{i+1}(\sigma^{d})$
\EndFor
\State append list of pairs $(\alpha,g^-_{i+1}(\beta))$, where $\beta\in B_0$ 
\EndFor
\EndFor
\end{algorithmic}
\end{algorithm}

Births and deaths are now defined in the same way  as in the case of a single cell decomposition.

\begin{defn}\label{borndied2} Suppose $\alpha$ is critical for $V_i$. If $\alpha$ is not strongly connected to any critical cell for $V_{i+1}$, we say $\alpha$ {\em dies} at $i$. If $\alpha$ for $V_i$ is not strongly connected
to any critical cell for $V_{i-1}$ we say $\alpha$ is {\em born} at $i$. If $\alpha$ is strongly connected to exactly one critical cell $\beta$ for $V_{i+1}$
and $\beta$ is not strongly connected to any other
critical cell for $V_i$ we  say $\alpha$ {\em moves} to $\beta$.
\end{defn}

As before, we may then draw diagrams illustrating these connectivity data. The algorithm for constructing these is exactly the same (Algorithm \ref{birthdeathalg2} below).

\begin{algorithm}
\caption{Algorithm Birth-death \#2}
\label{birthdeathalg2}
\begin{algorithmic}[1]
\For{$d=0,\ldots,\dim M$}
\For{$i=0,\ldots r$}
\ForAll{$\alpha^d$ critical in $V_i$} 
\State add a node (denoted by $\alpha$) in layer $i$ 
\If{$i>0$}
\ForAll{$\alpha^d$ critical in $V_i$}
\State add an edge from the node $\alpha$ to each node in $V_{i-1}$, strongly connected to $\alpha$ according to Algorithm \ref{fc2} 
\EndFor
\EndIf
\EndFor
\EndFor
\EndFor
\end{algorithmic}
\end{algorithm}

\begin{figure}
\centerline{\includegraphics[width=4.5in]{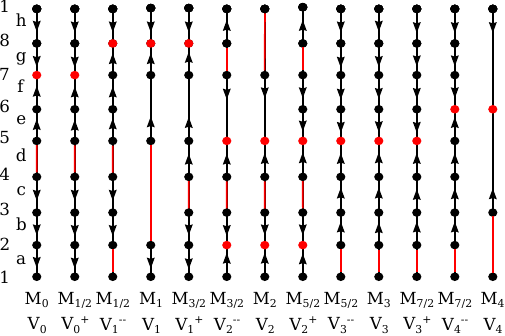}}
\caption{\label{alg2fig1} A family of discrete gradient fields on the circle (different triangulations).}
\end{figure}

\begin{figure}
\centering
\begin{subfigure}{.5\textwidth}
  \centering
  \includegraphics[width=.9\linewidth]{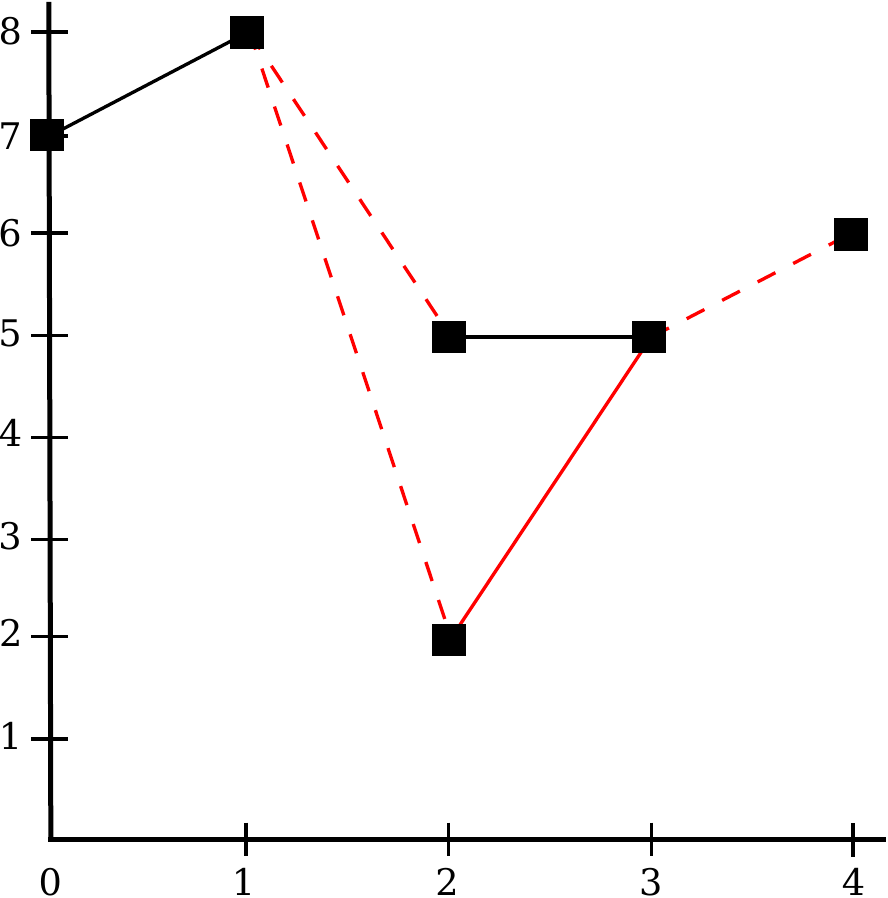}
  \caption{$0$-cells}
\end{subfigure}%
\begin{subfigure}{.5\textwidth}
  \centering
  \includegraphics[width=.9\linewidth]{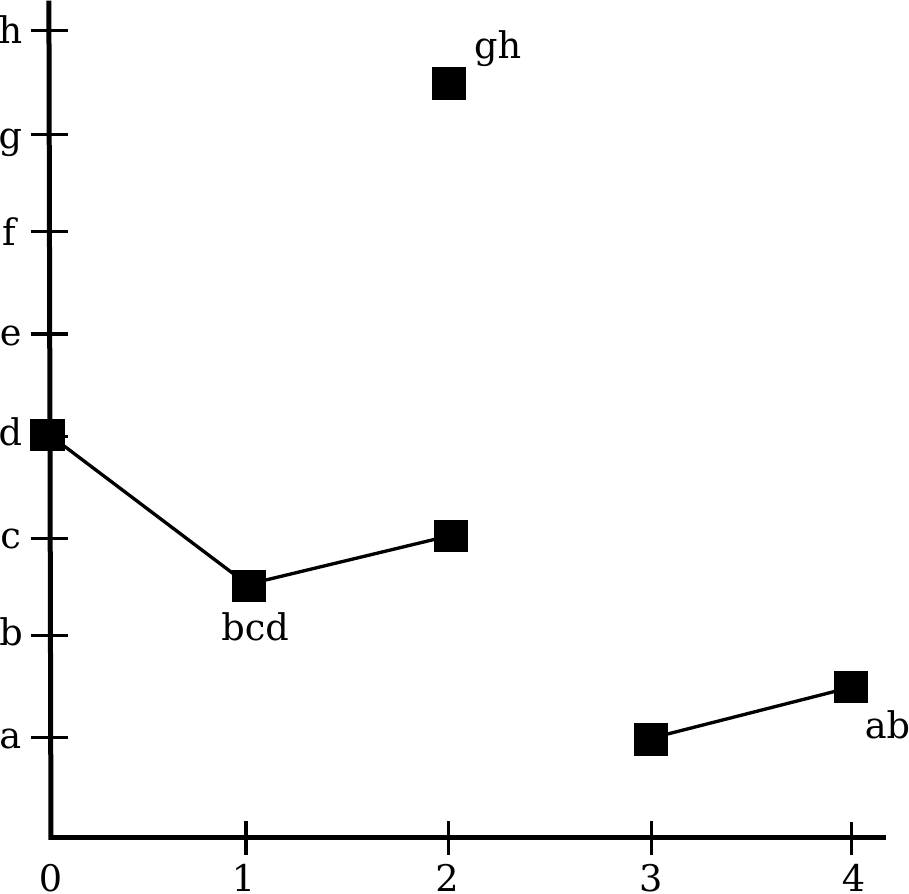}
  \caption{$1$-cells}
\end{subfigure}
\caption{\label{alg2fig2} The connectivity diagrams associated with the example in Figure \ref{alg2fig1}.}

\end{figure}

In Figure \ref{alg2fig1} we present an example of a family of discrete gradients on the circle, where each circle has a different triangulation.  We have indicated the refinements $M_{(2i\pm 1)/2}$ and the vector fields $V_{i\pm 1}^\pm$.  Figure \ref{alg2fig2} shows the connectivity diagrams for the $0$- and $1$-cells.  Strong connections are indicated by black lines, forward connections in adjacent slices are shown in red, and backward connections in adjacent slices are shown with red dashed lines. In this example we chose to show backward connections to illustrate how it is possible for a birth to occur, even though the new critical cell is connected to a critical cell in a previous slice. Our definition of birth requires a {\em strong connection}, the case of vertex 6 in slice $4$ shows this one-way connection. The births are: vertex $2$ and vertex $5$ in slice $2$ and vertex $6$ in slice $4$; edge $gh$ in slice $2$ and edge $a$ in slice $3$.  The deaths are: vertex $8$ in slice $1$, vertex $2$ in slice $2$, and vertex $5$ in slice $3$; edge $c$ in slice $2$ and edge $gh$ in slice $2$.

\section{Applications}\label{apps}

There are numerous applications of discrete Morse theory to data.  Here is a (rather vague) description of one possible application of the algorithms developed in this work. Take a network of sensors measuring some physical quantity in some domain (for example the temperature distribution in a gas or liquid, the presence of some chemical, light intensity, etc.). The sensors communicate with other sensors that are within their connectivity range, and the data is collected and processed by a central unit. The usual approach to sensor networks in topological data analysis is to model the domain of the network as a simplicial complex with the help of some topological reconstruction algorithm (for example the Vietoris--Rips or \v Cech complex). There is extensive literature on using this approach to the problem of domain coverage in sensor networks (\citep{GdS1},\citep{GdS2}). By extending the values of the function given on the vertices to a discrete Morse function (for example using the algorithm in \citep{kkm}), the overall distribution of the measured quantity can be analyzed. Algorithm \ref{birthdeathalg1} enables tracking significant features of the distribution for stationary sensor networks, while Algorithm \ref{birthdeathalg2} deals with time dependent sensor networks. 

In this section we describe an application of the birth-death algorithms to images. Birth-death diagrams produced by these algorithms can be used for tracing features through a sequence of images, as described below. A 2D (or 3D) grayscale digital image is given by a function $f :D\to \mathbf{Z}$ where $D$ is a rectangular subset of the lattice $\mathbf{Z}^2$ (or $\mathbf{Z}^3$ in case of 3D images) of size corresponding to the resolution of the image (i.e., the number of pixels or voxels) which associates to each point on the lattice the color intensity at the correspoding pixel (or voxel). For our purposes we model the image as a cubical complex with a vertex for each pixel and with pixels corresponding to neighbouring points on the lattice connected by edges and cubes of higher dimension. Figure \ref{rasters} (top row) shows two $9\times9$ rasters together with the corresponding cell complexes. The grayscale values, given in the vertices of the complex, can be extended to a discrete gradient vector field using one of the existing algorithms \citep{kkm}, \citep{rws}. 


In our first experiment the algorithm was tested on a sequences of images obtained from a CT scan of the head where the
algorithm is used to trace minima of the density function which
correspond to cavities and tunnels.  The cell complex underlying each slice is a $2$-dimensional cubical complex of voxelized data. The resolution of the images is $103\times 103$.

The number of local minima of the density function in one image
ranges from 200 to 300, and most of these correspond to small
fluctuations in the density and to noise.  In order to eliminate
these, as well as local minima with higher function values which
correspond to softer tissue, we use persistence.  Persistence is
incorporated in the algorithm for generating a discrete vector
field and enables canceling pairs of critical cells that are
connected by just one gradient path and have values that differ by
less than the chosen persistence level $p$.  By setting the
persistence level to a suitable value, pairs of critical points
which correspond to a small fluctuation in the density function,
as well as local minima which represent softer tissue are
cancelled, while local minima with low values of the density
function which correspond to the cavities are preserved.

\begin{figure}[h!]
\centerline{\includegraphics[height=9cm]{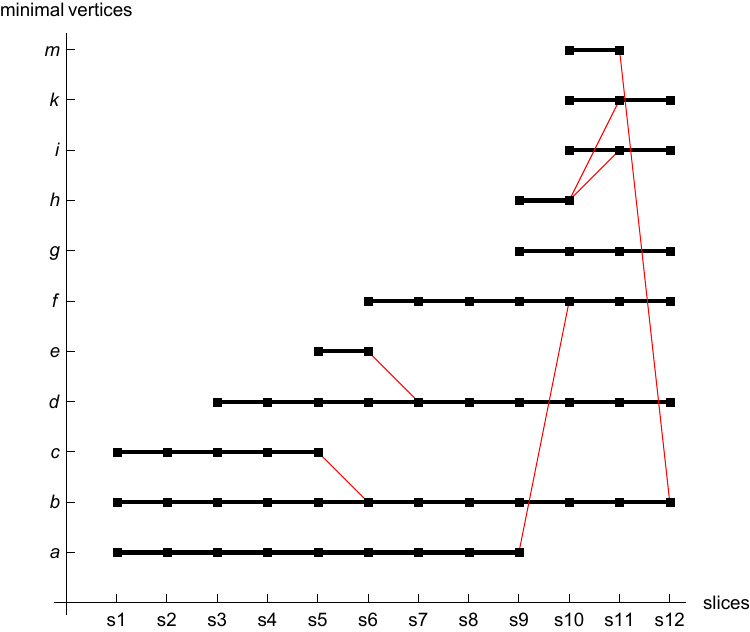}}
\caption{\label{dgm}Births and deaths of minima in a sequence of
head slices.}
\end{figure}

Figure \ref{dgm} shows the birth and death diagram of the minimal
vertices in a sequence of 12 images obtained by running the
algorithm with persistence set to $p=30$ (where the values range
from 0 to 255). Strong connections correspond to black lines, and connections that appear in one direction (in this case forward) but not in both correspond to red lines.

\begin{figure}[h!]
\centerline{\includegraphics[height=4cm]{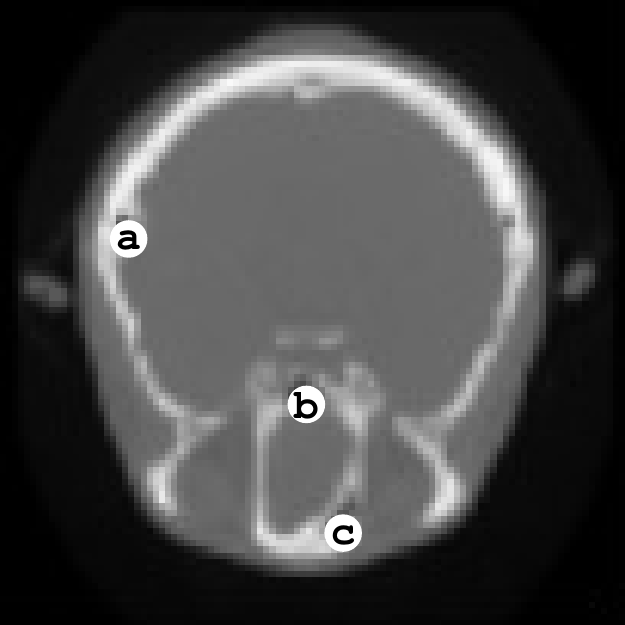}
\includegraphics[height=4cm]{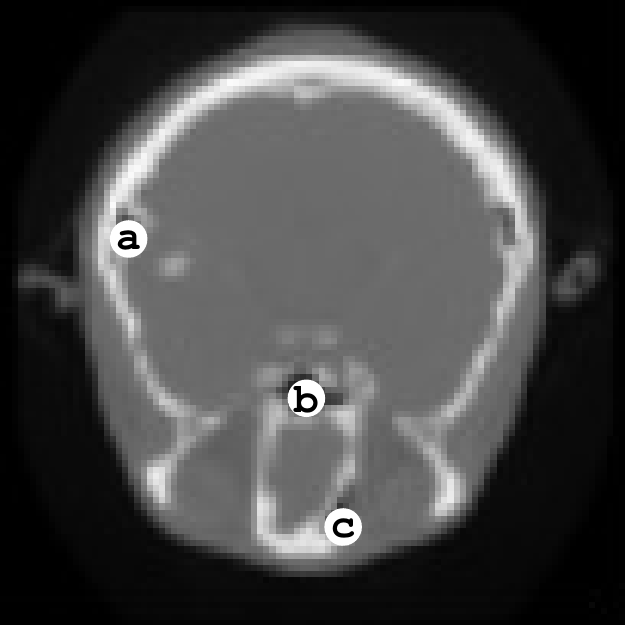}
\includegraphics[height=4cm]{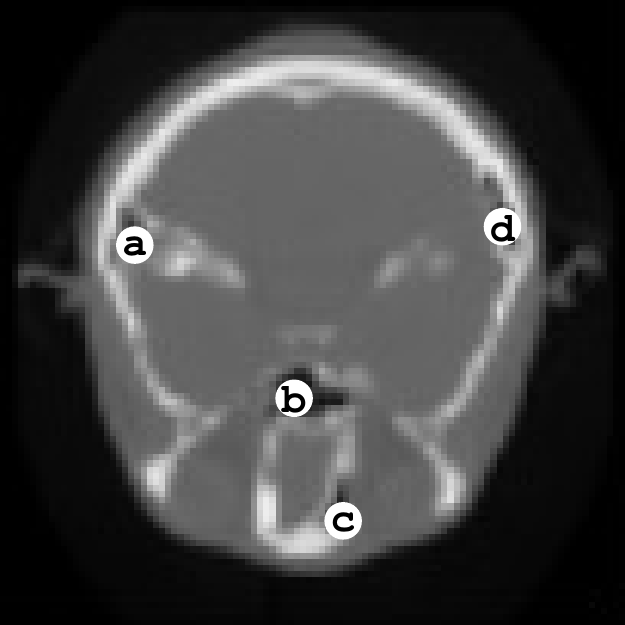}
\includegraphics[height=4cm]{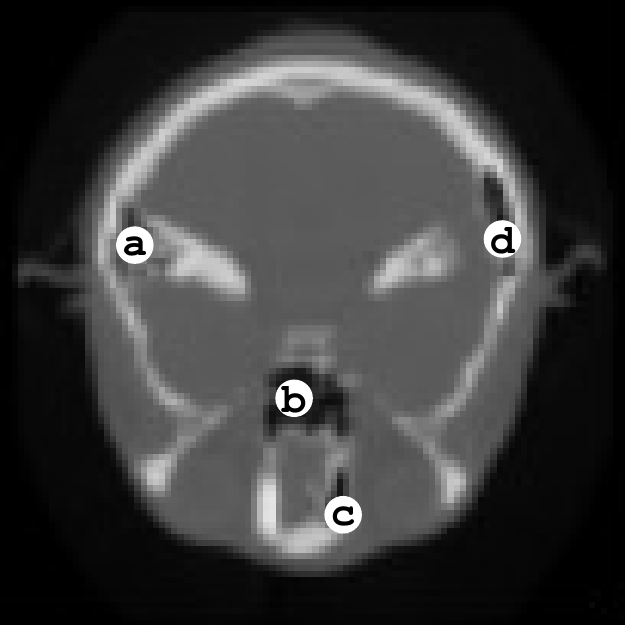}}
\vskip 3mm \centerline{\includegraphics[height=4cm]{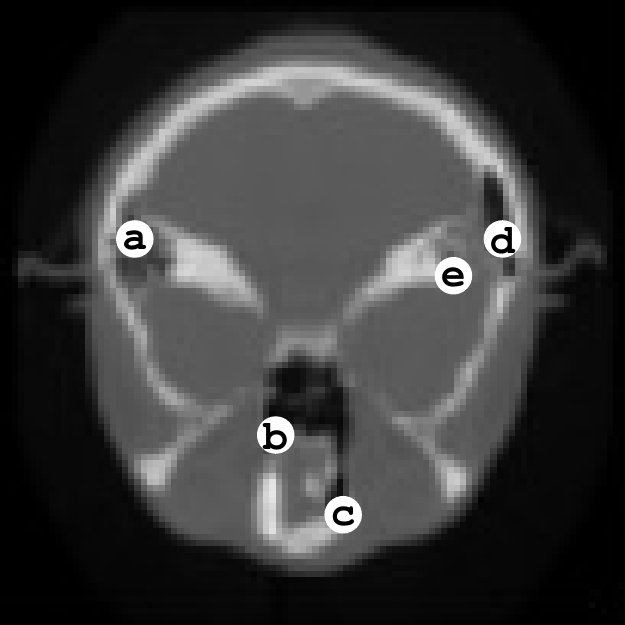}
\includegraphics[height=4cm]{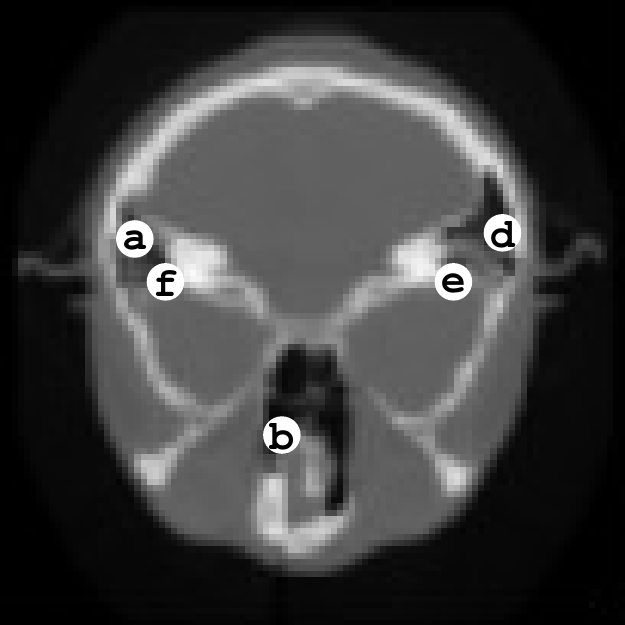}
\includegraphics[height=4cm]{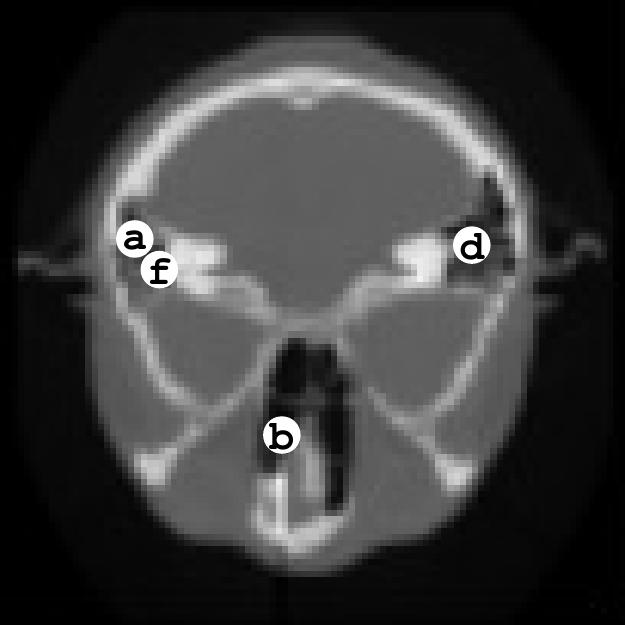}
\includegraphics[height=4cm]{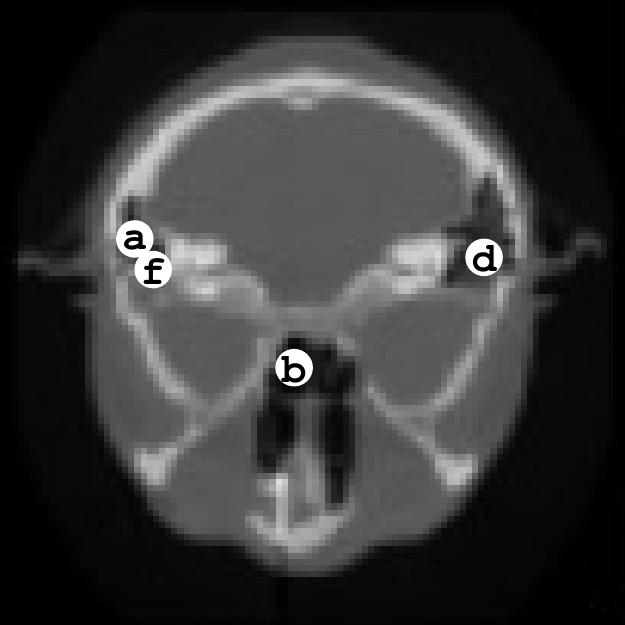}}
\vskip 3mm \centerline{\includegraphics[height=4cm]{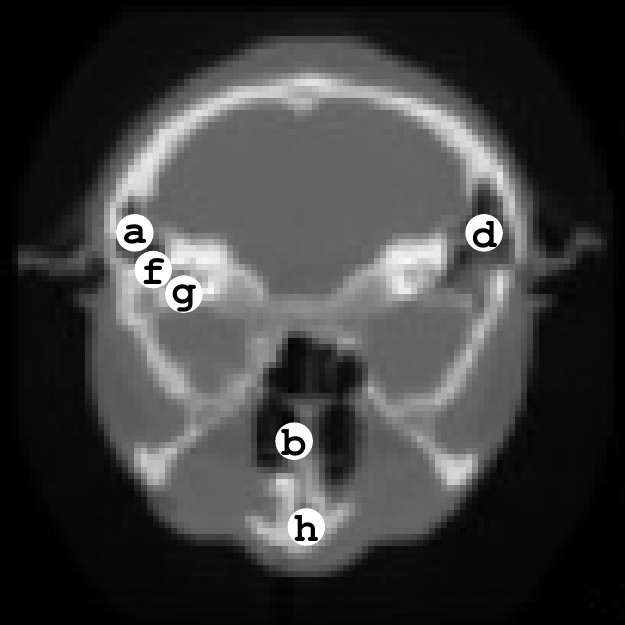}
\includegraphics[height=4cm]{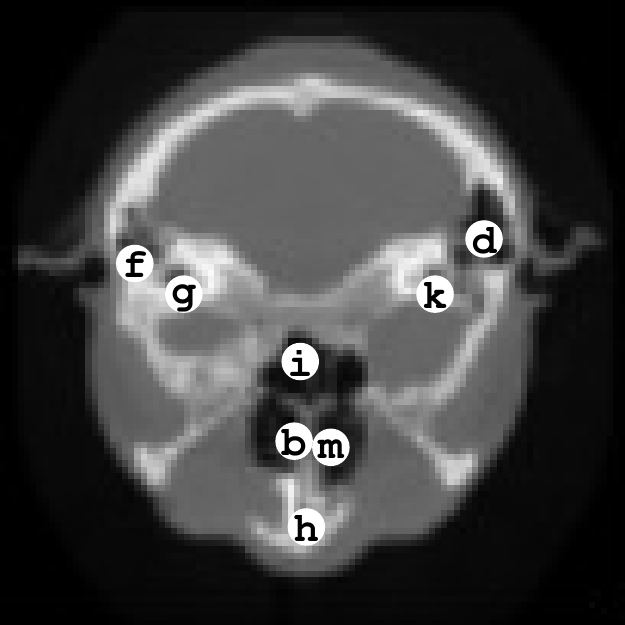}
\includegraphics[height=4cm]{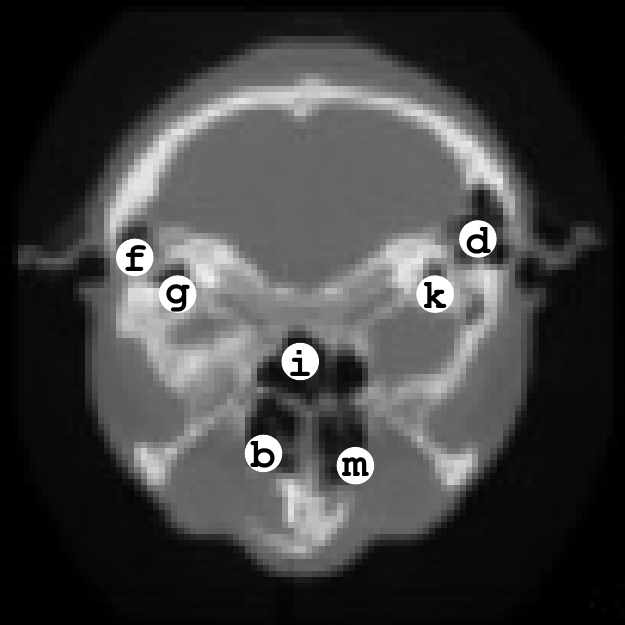}
\includegraphics[height=4cm]{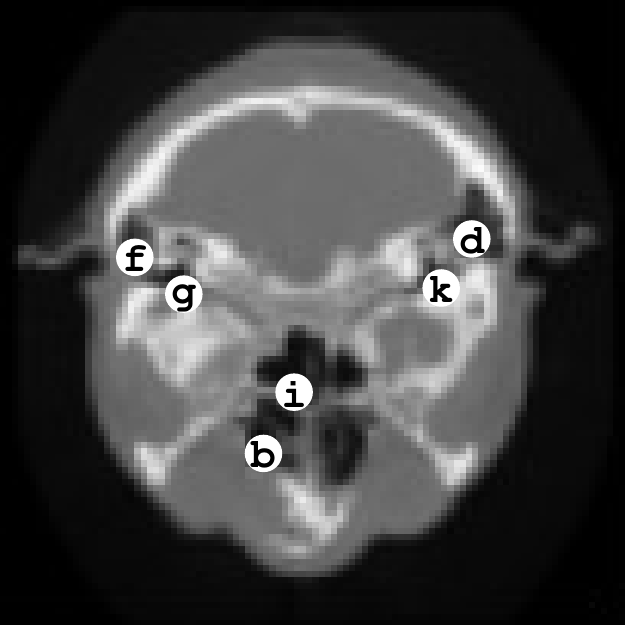}}
\caption{\label{slices}Tracing cavities through head CT scans.}
\end{figure}

\begin{figure}[h!]
\centerline{\includegraphics[width=4in]{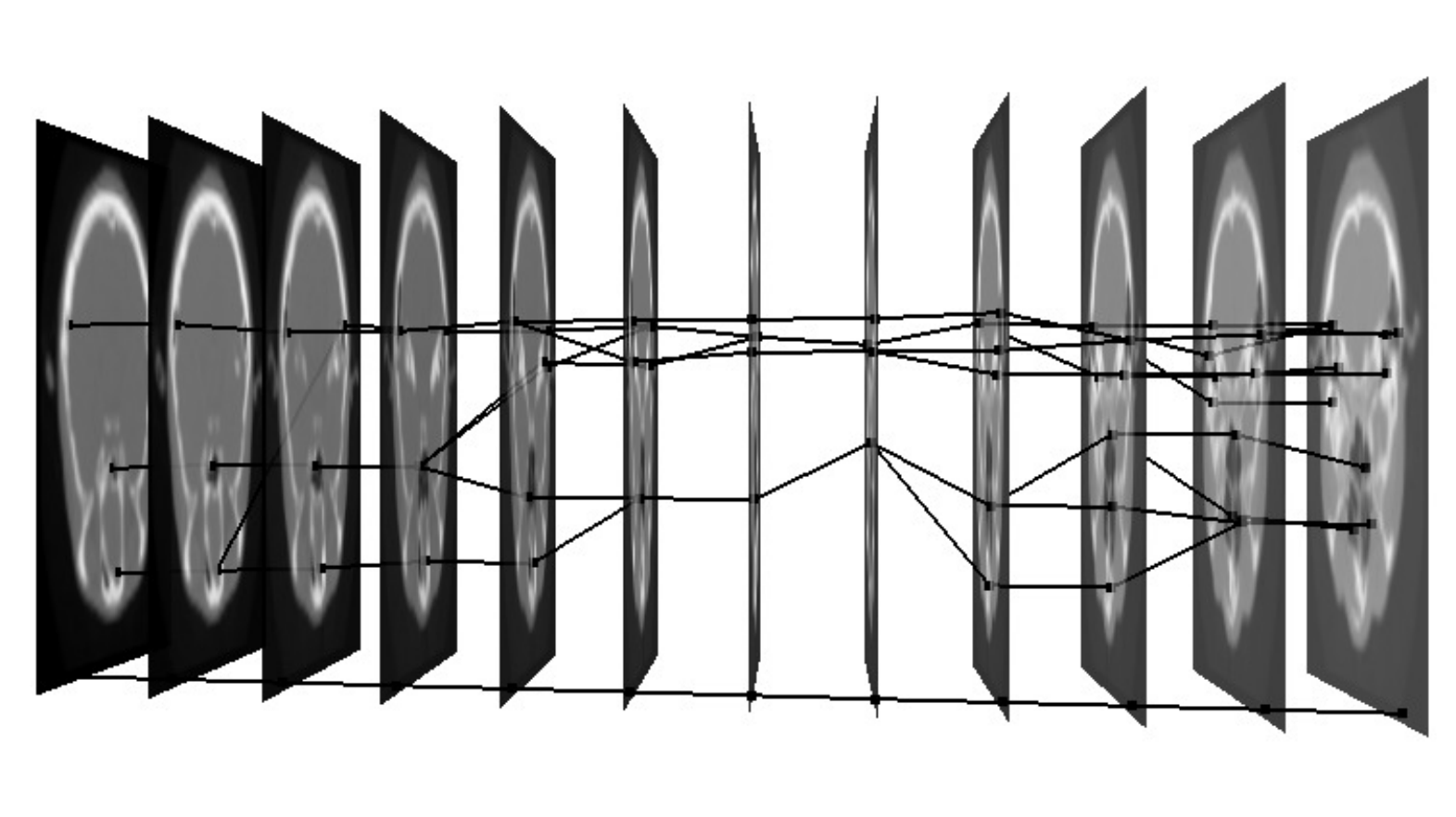}}
\caption{\label{tunnels}Simulated cavities in the head.}
\end{figure}

Figure \ref{slices} shows the sequence of images in which the
cavities connected by strong connections are traced. Points which
are strongly connected are labeled by the same letter. Altogether
12 strong connections have been found. As can be seen, each strong
connection represents a cavity, and most cavities are detected. Figure \ref{tunnels} shows the connections among the local minima, simulating the cavities themselves.

The second experiment deals with a series of scintigrams
obtained in renal scintigraphy. Images of the abdomen are taken
with a gamma camera at fixed time intervals after a radioactive
agent has been injected. The purpose is to follow the flow of the
agent through the kidneys in order to evaluate their function. The
algorithm was used to trace points of maximal concentration of the
agent which correspond to locations within the kidneys and the
bladder.  Again, the underlying cell complex in each slice is a $2$-dimensional cubical complex of voxelized data.

\begin{figure}[h!]
\centerline{\includegraphics[height=4cm]{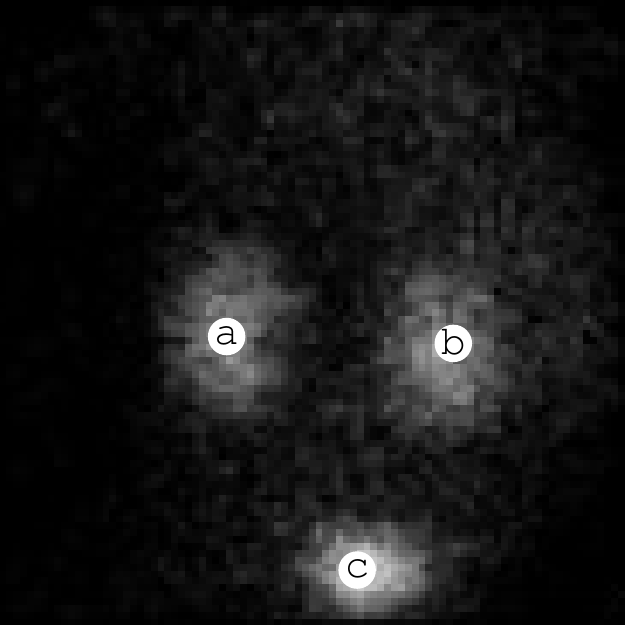}
\includegraphics[height=4cm]{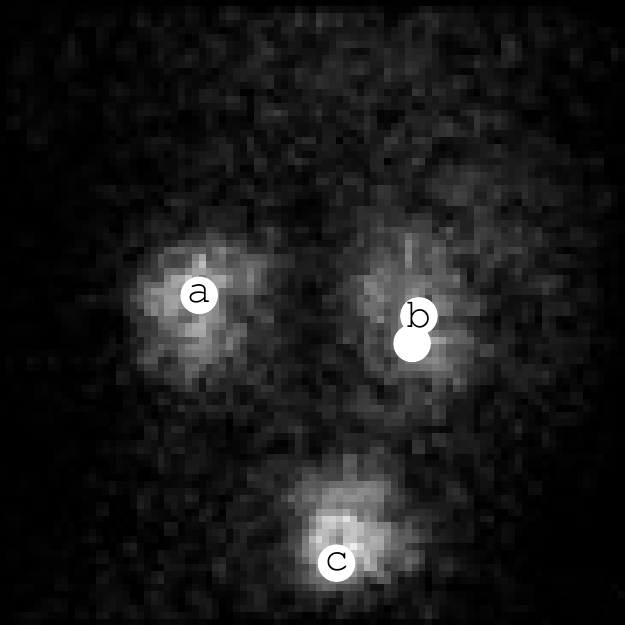}
\includegraphics[height=4cm]{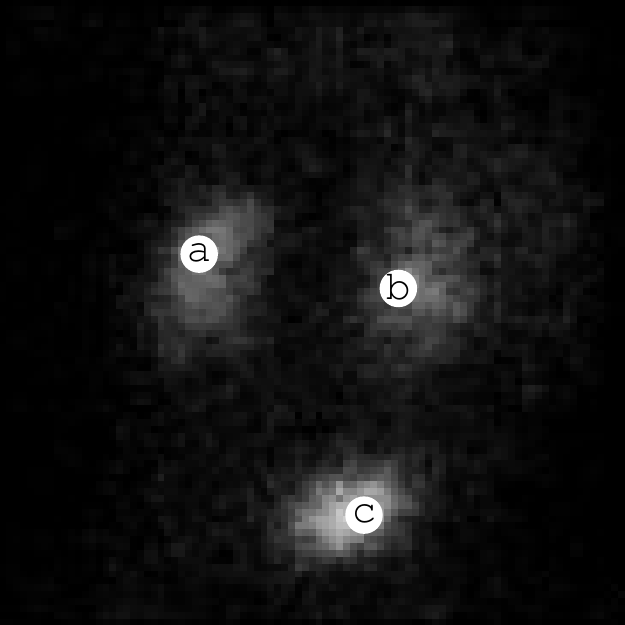}
\includegraphics[height=4cm]{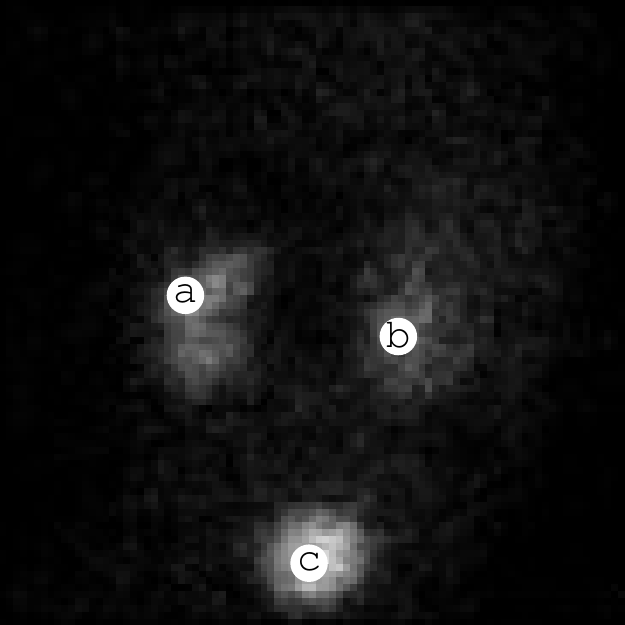}}
\vskip 3mm \centerline{\includegraphics[height=4cm]{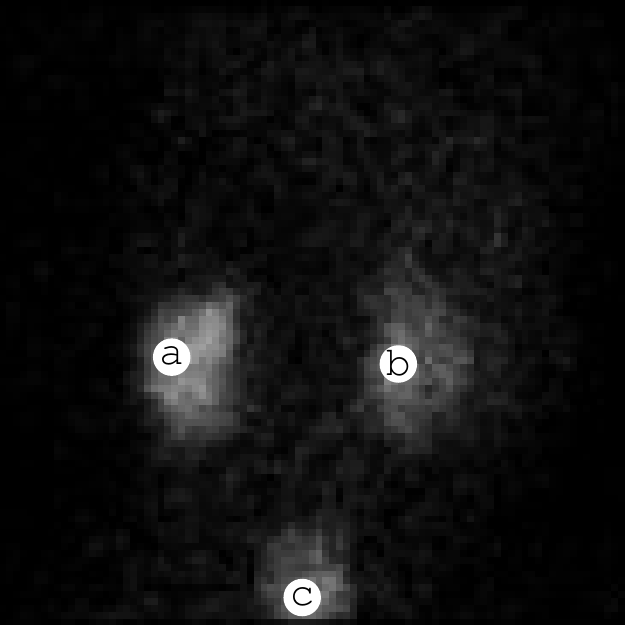}
\includegraphics[height=4cm]{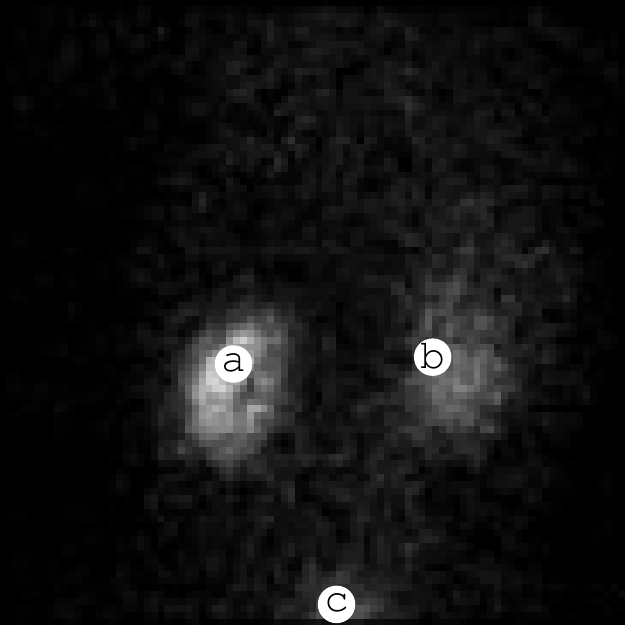}
\includegraphics[height=4cm]{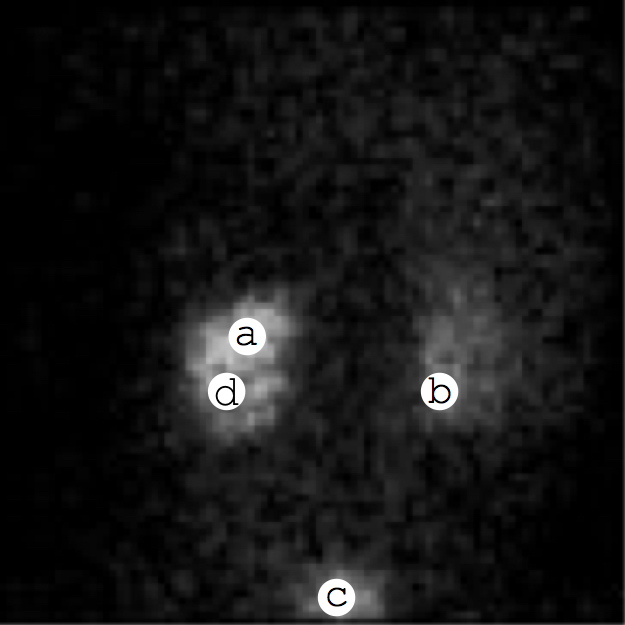}
\includegraphics[height=4cm]{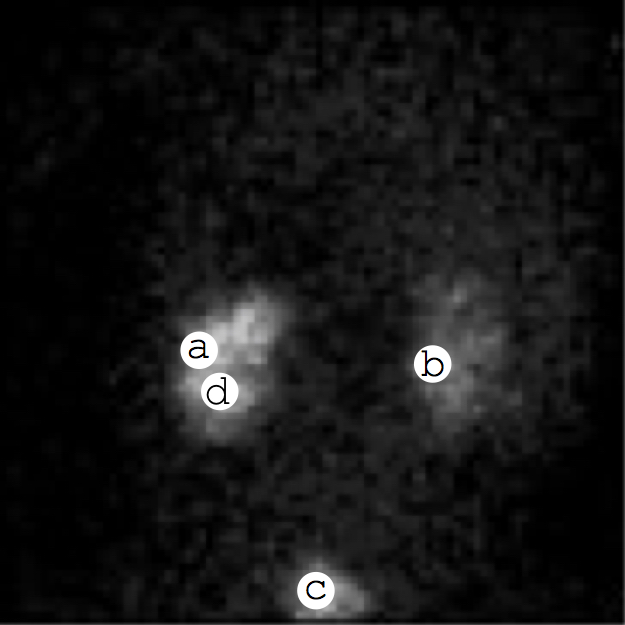}}
\caption{\label{kidneys}Tracing the maxima of radioactivity in
scintigrams.}
\end{figure}

\begin{figure}[h!]
\centerline{\includegraphics[height=7cm]{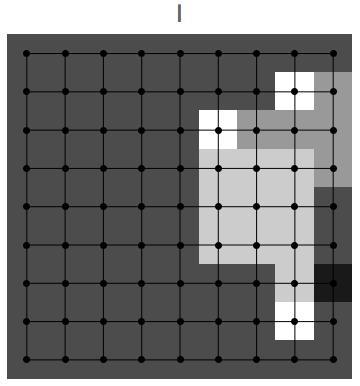}\quad\quad\quad\includegraphics[height=7cm]{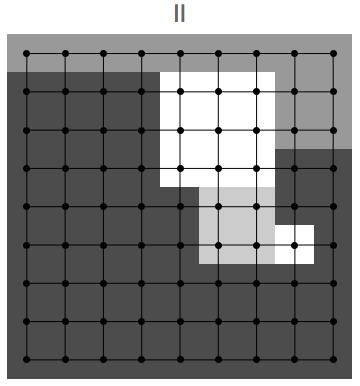}}
\centerline{\includegraphics[height=7cm]{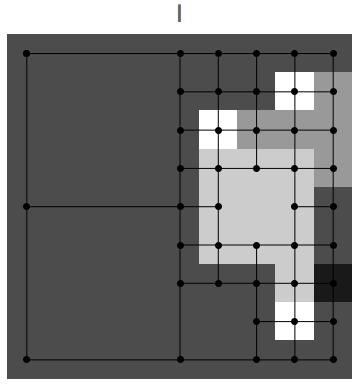}\quad\quad\quad\includegraphics[height=7cm]{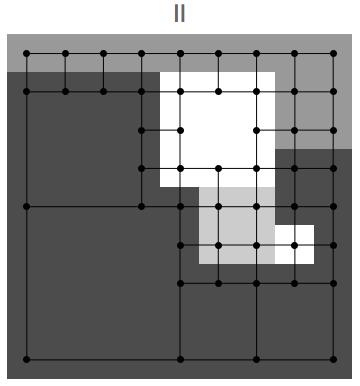}}
\caption{The top row shows two $9\times 9$ rasters with the corresponding initial cubical cell complex, and the bottom row with the reduced cell complexes.\label{rasters}}
\end{figure}

In contrast to CT scans, where noise is present in the form of a
relatively small error in the densities, the gray-level function
on a scintigram is heavily polluted by a quantum-type noise, where
the size of the error in a particular point can be relatively
large. Eliminating noise represents a major problem in the
analysis of this type of data set, and a standard approach in
image analysis is to use a smoothing convolution filter before
processing the images. Persistence canceling provides a different
way of avoiding this problem.  Figure \ref{kidneys} shows how the
algorithm traces the maxima of the grayscale function on
unfiltered images. The persistence level has been set to 50 (on a
scale from 0 to 255).

Our experiment shows that running the algorithm on images smoothed
by a Gaussian filter and setting the persistence to a lower value
of $p=10$ (an even lower value, for example $p=0$, gives too many
local maxima which correspond to outliers despite the smoothing,
so it is not good for tracing the locations) produces a similar
pattern as the one in Figure \ref{kidneys}, but typically with
only one peak within one kidney. Multiple peaks like those
appearing in the last two images are also of interest,
though, since they give insight into the activity within the
kidney. So prefiltering images seems to cause loss of useful
information and using unpreprocessed images and higher persistence
is a better choice.

Our Algorithm \ref{birthdeathalg1} has also been used in a machine learning algorithm from a robotic domain where a robot observing a red and a blue ball with an on-board camera is taught the concept of occlusion from examples of image sequences \citep{zabkar}.

\medskip

A major problem in using discrete Morse theory on images is that images often contain large areas which are more or less of the same color, that is, the  grayscale value is almost constant. Morse theory (discrete as well as smooth) is not well-suited for dealing with such functions and as a result algorithms for generating a discrete vector field typically produce a large number of unnecessary critical cells, mostly due to noise. This problem can be eliminated by a preprocessing step, where cubes of maximal dimension that cover regions with uniform grayscale values are merged into bigger cubes, reducing the complexity of the cell complex encoding the image.  Recall that pixels/voxels in the image are represented by vertices of the cubical complex with associated grayscale values. Specifically, the vertices (i.e. pixels) are  processed one by one and the complex is updated at each step in the following way. For each vertex the largest cube containing this vertex with grayscale values in all vertices differing by less than the preferred threshold is constructed. This large cube is then added to the complex, substituting all smaller cubes contained in it, and the boundary information is updated accordingly. The resulting complex consists of cubes of different sizes. 

In a sequence of images preprocessed in this way, critical cells can be connected using Algorithm 4. We demonstrate this on the simple example of the two $9\times 9$ rasters shown on Figure \ref{rasters}. In the top row the two rasters, denoted by I and II are shown together with the initial cubical complex  encoding a $9\times9$ raster, and in the botton row the reduced complexes are shown.


\begin{figure}
\centerline{\includegraphics[height=7cm]{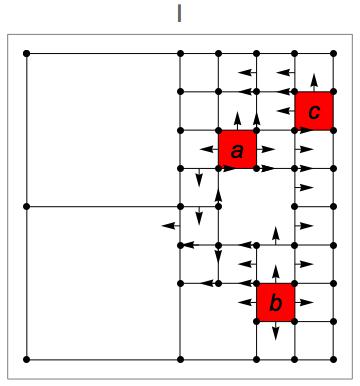}\quad\quad\quad\includegraphics[height=7cm]{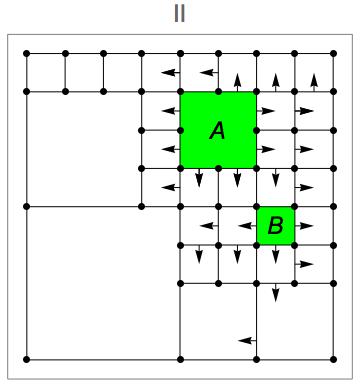}}
\centerline{\includegraphics[height=7cm]{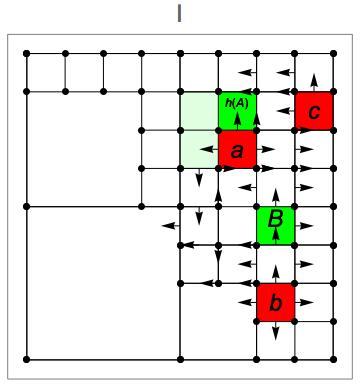}\quad\quad\quad\includegraphics[height=7cm]{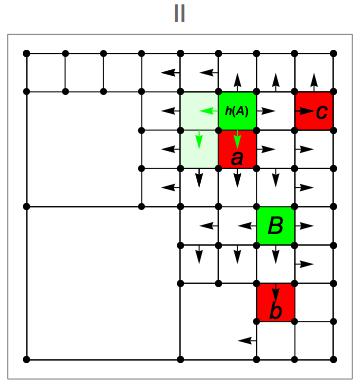}}
\caption{Top row: partial discrete vector fields on the complexes of rasters I and II from Figure \ref{rasters}. Bottom row, left: forward connections from I to II are $a\mapsto A$, $b\mapsto B$; right: backward connections are $A\mapsto a$, $A\mapsto c$ and $B\mapsto b$.\label{dvfs}}
\end{figure}

Figure \ref{dvfs} shows in the top row (a part of) the discrete vector fields obtained by extending the grayscale values given on vertices of the complexes on Figure \ref{rasters}. On raster I there are three critical 2-cells labeled by $a,b$ and $c$ and colored red, while on raster II there are two critical 2-cells labeled $A$ and $B$ and colored green. In the second row the common refinement of the cellulations is shown, with the refined discrete vector field of I on the left and the refined discrete vector field of II on the right. On both images, critical 2-cells from both slices are shown, where $h(A)$ denotes the refined critical cell corresponding to $A$ from Definition \ref{dvfrefinement}.  For forward connections from slice I to slice II, gradient paths on I starting in the red critical cells and leading to the green critical cells are followed. This gives us the connections $a\mapsto A$ (since there is a path starting in the boundary of $a$ leading to $h(A)$, and $b\mapsto B$. The critical cell $c$ of I is not forward connected to any green critical cell, so it dies at slice I. For backward connections, gradient paths on II should be followed, which give the connections $A\mapsto a$, $A\mapsto c$ and $B\mapsto b$, so the critical cell $A$ of slice II splits backwards into the two critical cells $a$ and $c$. We therefore have strong connections between $a$ and $A$ and between $b$ and $B$.

\begin{figure}[h!]
\centerline{\includegraphics[height=5cm]{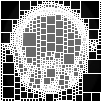}}
\caption{Reduced cell decomposition of a head CT scan. \label{cubedhead}}
\end{figure}

Note that the reduction in size of the complexes after merging depends on the type of the images, as well as on the heuristics used for merging cells. As an example, Figure \ref{cubedhead} shows the reduced decomposition of a CT head scan, where the number of cells is reduced by approximately 1/2. 

An extension of the algorithm presented in \citep{kkm} to polyhedral complexes was used for extending the grayscale values on the reduced complexes corresponding to the individual slices to discrete vector fields. For 2D images, the refinements of the discrete vector fields to common subdivisions have critical cells corresponding to the critical cells of the original discrete vector fields, and each path in the refined discrete vector field is covered by a path in the original discrete vector field. This implies that if threshold $0$ is used for merging cells in the original images, the resulting birth-death diagram on reduced complexes will be the same as the diagram obtained on the original images (Figure \ref{dgm}), while the computational complexity will be, in some cases, much smaller.  This illustrates the idea that even if we have the same cell decomposition for a space for multiple times, it may be advantageous to reduce the size of the decompositions first (via this merging technique, for example) and then apply Algorithms \ref{fc2} and \ref{birthdeathalg2}.

\section{Open questions and future work}\label{smooth}

It is natural to ask the following.

\begin{question} To what extent do the birth-death diagrams our algorithm
produces agree with those in the smooth case?
\end{question}

That is, if we have a family
of smooth maps $G:X\times I\to \zr$ which are generically Morse, and after
choosing a cell decomposition of $X$ and restricting the maps $G(-,t)$ to a finite
collection of slices we run our algorithm, how is the resulting birth-death
diagram related to that for $G$?  We will not answer this question completely.  Rather, we will provide evidence that suggests our method is a reasonable approximation.  To begin the study of this we first ask the following question.  Given a Morse function $F$ on a manifold $X$, does there
exist a cell decomposition of $X$ and a discrete Morse function on that
decomposition whose critical cells are in one-to-one correspondence with the
critical points of $F$?  This is answered by the following theorem of Gallais
\citep{gallais}.

\begin{thm}\label{onetoonecrit}{\em (c.f. \citep[Theorem 3.1]{gallais})} Let $X$ be
a smooth manifold and let $F:X\to\zr$ be a Morse function with gradient-like vector field $v$.  Assume that $F$ is generic in the sense that its ascending and descending manifolds are transverse.
Then there is a $C^1$-triangulation of $X$ and a \dgvf\ $V$
such that 

\begin{enumerate}
\item the critical cells of  $V$ are in one-to-one correspondence with the
critical points of $F$;

\item for each pair of critical cells $\sigma_p$ and $\sigma_q$ of $V$ such that $\dim(\sigma_p) = \dim(\sigma_q)+1$, $V$-paths from the boundary of $\sigma_p$ to $\sigma_q$ are in bijection with integral curves of $v$ up to a renormalization connecting $q$ to $p$.
\end{enumerate}
Moreover, the critical points lie in the interiors of
the corresponding critical cells and the closures of critical cells are all disjoint.  
\end{thm}

\begin{proof}
All except the disjointness of the critical cells is a consequence of
Theorem 3.1 of \citep{gallais}.  But the critical cells can be made disjoint by taking barycentric subdivisions
and applying Proposition 6.24 of \citep{knudson}.
\end{proof}

\begin{Rm} There is a small gap in the proof of Gallais's theorem, but this has recently been fixed by Benedetti \citep{benedetti}, who also proved a stronger version of the result.
\end{Rm}

Now suppose we have a family of such maps and we choose a finite collection of slices and the Morse function $F_i$ on each.  Theorem \ref{onetoonecrit} gives us a triangulation $M_i$ of each slice and a discrete gradient $V_i$ whose critical cells are in one-to-one correspondence with the critical points of $F_i$.  We would then employ Algorithm \ref{birthdeathalg2}.  This requires choosing common refinements of adjacent slices and finding refinements of the $V_i$.  We want to do this without generating additional critical cells; that is, we seek refinements $V_i^\pm$ with the same numbers of critical cells as $V_i$ and $V_{i-1}$, respectively.  Moreover, we would like the critical cells of the refinements to be contained in the critical cells of $V_i$ and $V_{i-1}$.  This is in fact possible; Proposition 6.24 of \citep{knudson} shows how.  
  Since the critical cells of the various discrete vector fields correspond to the critical points of the Morse functions and the $V_i$-paths correspond to integral curves of gradient-like vector fields for the functions, the birth-death diagrams should approximate the smooth bifurcation diagrams, provided the sampling is sufficiently fine.  

Recall that Cerf's Theorem on $1$-parameter families of smooth maps \citep{cerf} allows us to assume the following.  A $1$-parameter family of smooth functions on a smooth compact manifold $X$ is uniformly close to a family $G:X\times I \to\zr$ satisfying the following properties:
\begin{enumerate}
\item The function $G(-,t)$ is Morse at all but finitely many values of $t$;
\item If the function $G(-,t)$ is not Morse, then it has a single nondegenerate critical point $p$, and near that point the family $G(-,t)$ is conjugate to the family $f_t(x_1,\dots ,x_n) = G(p,t) + x_1^3+ \varepsilon_1 tx_1 + \varepsilon_2 x_2^2 +\cdots + \varepsilon_n x_n^2$, with $\varepsilon _i=\pm 1$.
\end{enumerate}
In the case $\varepsilon_1 = -1$, there is a birth, and in the case $\varepsilon_1=1$, there is a death.

If there is to be any chance for the diagrams generated by Algorithm \ref{birthdeathalg1} and Algorithm \ref{birthdeathalg2} to mimic with those built from a smooth family, we will be forced to make some density assumptions about the sampling.  Suppose we have a family $G:X\times I \to\zr$ of the type guaranteed by Cerf's Theorem.  We want to assume that the slices are close enough to ensure that the critical points of the Morse functions do not move very much.  This is certainly possible, but then it is not clear how to translate this idea to the discrete case.  What is really needed to formalize this is the following.

\begin{question} Is there a useful theory of ``closeness" for discrete vector fields?
\end{question}

Developing this would certainly be interesting and useful, but it is beyond the scope of this paper.  We therefore content ourselves with heuristics and note that the applications in Section \ref{apps} show the utility of our algorithms.

\begin{question}  Are there conditions that guarantee that a critical cell in one slice is strongly connected to at most one critical cell in the next slice?
\end{question}

Such a result would imply that the discrete diagrams generated by our algorithms are a good approximation to the smooth case.  It may also have independent utility for data analysis problems.

\subsection*{Acknowledgements} We are grateful to Gregor Jer\v{s}e for his help with the examples in Section \ref{apps}.  We also thank a pair of anonymous referees whose comments substantially improved the readability of the paper.

\end{document}